\documentclass[12pt]{article}
\usepackage{epsfig}
\usepackage{latexsym}
\usepackage{amsmath}
\usepackage{amssymb}
\usepackage[english]{babel}
\usepackage{graphicx}
\DeclareGraphicsExtensions{.eps,.bmp}
\usepackage[utf8]{inputenc}
\usepackage{float}
\usepackage{epsfig}
\usepackage{eepic}
\usepackage{pict2e}
\usepackage{xcolor}
\DeclareMathOperator*{\Ls}{\mathrm{Ls}}
\DeclareMathOperator*{\Linf}{\mathrm{Li}}
\DeclareMathOperator*{\Lim}{\mathrm{Lim}}

\usepackage{ifpdf}

\newtheorem{thm}{Theorem}[subsection]
\newtheorem{cor}[thm]{Corollary}
\newtheorem{lem}[thm]{Lemma}
\newtheorem{prop}[thm]{Proposition}
\newtheorem{defn}[thm]{Definition}

\newcommand{\un}{1\!\!1}

\newcommand{\Real}{\mathbb{R}}

\numberwithin{equation}{section}

\numberwithin{table}{section}
\numberwithin{figure}{section}

\newfont{\sfl}{cmssi12}

\usepackage[left=2.5cm,top=2.5cm,right=3.2cm,nohead,foot=2.5cm]{geometry}

\begin{document}

\title{ Kolm-Pollack Form, Translation Homotheticity and Tropical Limit of Production Technologies\thanks{The authors greatly acknowledge the French Research National Agency for the LAWBOT Project: ANR-20-CE38-0013.}}

\author{Walter Briec\thanks {University of Perpignan, Department of Economics 52 avenue Villeneuve,
66000 Perpignan, France.}  \\ \textsc{lamps},
Universit\'{e} de Perpignan  \and St\'{e}phane
Mussard\thanks{$^{{}}$\textsc{Chrome} Univ. Nîmes, Rue du Dr Georges Salan, 30000 Nîmes, - e-mail: stephane.mussard@unimes.fr, Research fellow \textsc{Um6p} Polythecnic Rabat, \textsc{Liser} Luxembourg.} \\ \textsc{Univ. Nîmes Chrome} , \textsc{Um6p}
\and Paola Ravelojaona \\ ICN, 
Nancy}


\maketitle

\begin{abstract}
 In this paper,   we consider a new class of generalized Convex structure and we investigate their tropical limits. Some properties are pointing out such that translation homotheticity and others ones allowing to consider the case of discrete production sets that are related to some specific dual forms. Along this line a general class of mathematical programs are derived and it is shown that they can be computed using standard methods. The proposed approach allows to deal with efficiency measures (output oriented or input oriented) on continuous and discrete data. 
\end{abstract}

{\bf AMS:} 06D50, 32F17\\


\break

\section{Introduction}							
The standard non-parametric production model initiated by \cite{ccr78} and \cite{bcc84} helps to understand production problems involving quantities represented by real values. This is actually due to the convexity hypothesis and to the nature of the measures used, one of the best known being the Debreu-Farrell measure (see \cite{d51}, \cite{f57}). The efficiency measure is also based, like the convexity hypothesis, on the notion of divisibility. For example, for an input measure, it evaluates the maximum contraction of an input vector allowing the efficient frontier to be reached. In this context it is therefore difficult to apprehend the situation in which the variables which are considered are integers instead of being real or simply defined over the field of rational numbers.

However, there are alternative models to the classic convex models which relax the convexity hypothesis. One of them is based on the hypothesis of minimum extrapolation of the technology consisting in considering the smallest production set which contains all the production vectors observed and satisfying the hypothesis of free disposal, \textit{i.e.} an increase in inputs makes it possible to produce at least the same quantities of output. One of these models is known as FDH (``Free Disposable Hull, see \cite{dst84}) and it is suitable for taking into account integer variables. One of its limitations, however, is that it assumes quite rudimentary technology involving zero or infinite marginal productivities. Moreover, a radial measure involving divisibility of the benchmarks associated with the evaluation of efficiency may not have integer values. Perhaps more importantly, these technologies are not always very discriminating, as many firms can be efficient, especially when the number of inputs and outputs is large.

More recently, other types of non-convex production models have been proposed and are based on an algebraic structure involving a lattice and divisibility assumption. These models are called $\mathbb B$-convex (see \cite{bh04,bh08}, \cite{bl11}, \cite{abm17}) and are related to the topological limit of the CES-CET models (Constant-Elasticity-of-Substitution and Constant-Elasticity-of-Transformation, see \cite{fgn88}). However they imply divisibility like the standard DEA models \cite{bcc84, bm86}. 

In this paper we consider the tropical approaches proposed by \cite{abs17}. A Max-Plus approach was implicitly considered and it can be shown to be perfectly suited for the use of a new class of production models with homogeneous efficiency measures, the so-called  directional distance function (see \cite{bd80}, \cite{ccf96}). These production models involve a tropical algebraic structure in which the operation $max$ replaces addition, and $+$ replaces  multiplication. These algebraic structures were analyzed and popularized by \cite{ms}. They are based on the notion of dequantization. The constructed algebraic structure forms an idempotent semigroup, the operation $\max$ being idempotent but not symmetrizable. First, we extend this model to a class of Min-Plus models and establish relationships between these two types of models. We also show that these models can be seen as topological limits of a general class of production models that are associated with the Kolm-Pollack form \cite{K76} traditionally used in social choice theory. The dequantization proposed by Maslov is therefore based on an algebraic structure considered independently in economic inequality theory. After relating these different approaches, we show that these tropical models are quite suitable, combined with the directional distance function, to take into account integer variables. The only restriction to be imposed is to consider a direction associated with the unit vector in the space of inputs or in that of outputs. When the observed quantities are integers, the distance functions are also integers and the resulting benchmarks are also integers. Among the particularities of this class of technologies, they make it possible to take into account graph translation homothetic  in the graph. This holds true for Kolm-Pollack forms although they are fundamentally nonlinear.

The paper is organized as follows. Section \ref{DEA} sets the tools employed for the use of production models. Section \ref{Power} sets the mathematical frameworks about generalized convexity and the concept of dequantization (Kolm-Pollack form). Some general results about tropical limit sets are provided. Section \ref{Quantized} introduces input and output oriented distance functions for non-parametric models as well for DEA models, and duality. Section \ref{Dequantized} introduced the tropical limits of production technologies with closed forms of distance functions (output oriented or input oriented). Discrete production models and their related distance functions defined on integer values  are proposed. 

\section{The Non-Parametric Production Model}\label{DEA}

The mathematical tools underlying generalized convexities are applied to production models.
Subsections 1, 2 and 3 are devoted to the exposition of the basic concepts:
the production technology, the methods used to estimate the production frontier, and by the way, the technology set.

\subsection{The Background of the Production Model}

\textit{Notations.} Let $\Real_{+}^{d}$ be the non-negative $d$-dimensional Euclidean space. For $z,w\in \Real_{+}^{d}$, we denote $z\leq w$ if, and only if, $z_{i}\leq w_{i}$ for all $i\in[d]$ where $[d]=\big\{1,\ldots,d\big\}$. Let $\mathbb{N}$ (respectively $\mathbb Z$) be the set of non-negative integers (non-positive integers). Let $\mathbf{e}$ (or $\mathbf{exp}$) be the elementwise exponential function $\mathbf{e}^z = (e^{z_1},\ldots,e^{z_d})$, and $\mathbf{ln}$ the elementwise logarithm function $\mathbf{ln}(z) = (\ln(z_1),\ldots,\ln(z_d))$. The vector $\un_d$ stands for the $d$-dimensional vector of ones. For all $m,n\in \mathbb{N}$, such that $d=m+n$, a production technology transforms
inputs $x=(x_{1},\ldots,x_{m})\in \mathbb R_+^m$ into outputs $y=(y_{1},\ldots,y_{n})\in \mathbb R_+^n$. The set $T\subset \mathbb{R}_{+}^{m+n}$ of all input-output
vectors that are feasible is called the production set.
It is defined as follows:
\begin{equation}\notag
T=\left\{ (x ,y )\in \mathbb{R}_{+}^{m+n}
: x \hbox{ can produce } y \right\}
\end{equation}
\noindent $T$ can also be characterized
by an input correspondence $L: \Real_+^n\longrightarrow 2^{\Real_+^m}$ and an output
correspondence $P:\Real_+^m\longrightarrow 2^{\Real_+^n}$ respectively defined by,
$$L(y)=
\left\{ x \in \mathbb{R}^{m}
: (x,y) \in T \right\} \text{ and }P(x )=
\left\{y\in \mathbb{R}^{n}
: (x,y) \in T  \right\}$$
The production set  $T$ can be identified with its graph, that is:
$$T={\displaystyle\big\{(x, y)\in \Real_+^m\times \Real_+^n: x\in L(y)\big\}}={\displaystyle\big\{(x, y)\in \Real_+^m\times \Real_+^n: y\in P(x)\big\}}$$ 
The inverse of  $P$  is the input correspondence   $L$   defined by
$x\in L(y)$ if and only if $y\in P(x)$. The sets $P(x)$ are
the values of $P$ while the sets $L( y)$ are the
fibers of $P$.  The  image of a subset $A$ of $ \Real_+^m$  by $
P$ is the set $P(A) = \bigcup_{x\in A}P(x)$. Finally, let
us denote the cone of free disposability as:
\begin{equation}\notag
K=\mathbb{R}_{+}^{m}\times
(-\mathbb{R}_{+}^{n})
\end{equation}
This cone plays an central role to characterize the free disposal assumption (T3) defined below. The usual assumptions generally imposed on the
production technology \cite{s53} are the following.

\medskip

\noindent T1: $T$ is a closed set.

\noindent T2: $T$ is a bounded set, \emph{i.e} for any
$z\in T$, $(z-K)\cap T$ is bounded.

\noindent T3: $T$ is strongly disposable, \emph{i.e.}
$T=(T+K)\cap \mathbb{R}_{+}^{m+n}$.

\medskip

\noindent Assumptions T1-T3 define a convex technology with freely
disposable inputs and outputs. The following subsection presents a classical way to estimate the production
technology. Let us define as $\mathcal T$ the class of all production sets satisfying axioms T1-T3.

\subsection{Non-Parametric Convex and Non-Convex Technology}

\noindent Following the works initiated by \cite{f57}, \emph{et al.} \cite{ccr78} and \cite{bcc84}, the production set is traditionally defined by the
convex hull that contains all observations under a free
disposal assumption. Suppose that
$A=\{(x_1,y_1),\ldots,(x_\ell,y_\ell)\}\subset \Real_+^{m+n}$ is a
finite set of $\ell$ production vectors. Let $Co(A)$ denotes the convex hull of $A$. From \cite{bcc84}, the production set under an assumption of variable
returns to scale is defined by,
\begin{equation}\notag
T_{V}(A)=(Co(A)+K)\cap \Real_{+}^{m+n},
\end{equation}
or equivalently, for any given vector $\ {t}=(t_1,\ldots,t_\ell)$, by
\begin{equation}\notag
T_{V}(A)=\Big\{(x,y)\in \mathbb{R}_{+}^{m+n}: x \geq \sum_{k=1}^{\ell} t_{k}\,x_{k},
y\leq \sum_{k=1}^{\ell} t_{k}\,y_{k}, {t} \geq 0, \sum_{k=1}^\ell
t_{k}=1\Big\}
\end{equation}
This approach is the so-called DEA method (Data Envelopment Analysis)
that leads to an operational definition of the production set. This subset represents some kind of convex hull of the observed production vectors.
In line with \cite{ccr78}, under an assumption of constant returns
to scale, the production set can also be represented by
the smallest convex cone containing all the observed firms. In such a
case the constraint $\sum_{k=1}^{\ell}
t_{k}=1$ is dropped from the above model,  and therefore the production set becomes:
\begin{equation}\notag
T_{C}(A)=\Big\{(x,y)\in \mathbb{R}_{+}^{m+n}: x \geq \sum_{k=1}^{\ell} t_{k}\,x_{k},
y\leq \sum_{k=1}^{\ell} t_{k}\,y_{k}, {t} \geq 0 \Big\}
\end{equation}
Technical efficiency can be measured by introducing the usual concept of input distance function
and finding the closest point to any observed firms on the boundary of the production
set. Accordingly, the problem of efficiency measurement can be readily
solved by linear programming. Among the most usual measures
of technical efficiency, the Farrell efficiency measure (\cite{f57} and \cite{d51}) is essentially the inverse of Shephard's distance function \cite{s53}. The input Farrell efficiency measure is the map $E_{\mathrm{in}}: \Real_+^{m+n}\times \mathcal T\longrightarrow \Real_+\cup \{\infty\}$ defined as follows:

\begin{equation}\notag
E_{\mathrm{in}}(x,y, T)=\inf\Big\{\lambda \geq 0: \big({\lambda}  x,y\big)\in T\Big\}
\end{equation}

\noindent It measures the greatest contraction of an input
vector until to reach the isoquant of the input correspondence,
and can be computed by linear programming. In the output case, the output Farrell efficiency measure is the map $E_{\mathrm{out}}: \Real_+^{m+n} \times \mathcal T\longrightarrow \Real_+\cup \{\infty\}$ defined as:

\begin{equation}\notag
E_{\mathrm{out}}(x,y, T)=\sup\Big\{\theta \geq 0 : \big( x, {\theta} y\big)\in T\Big\}
\end{equation}

\noindent It is also possible to exogenously set inputs and outputs to measure efficiency \cite{bmo86}.
It is possible to provide a non-parametric estimation
that does not postulate the convexity of the technology. It is the FDH (Free Disposal Hull) approach developed by  \cite{dst84}. The FDH hull of a data set yields the following
non-parametric production set,
$$T_{F }(A)=(A+K)\cap \Real_{+}^{n+m}$$
The main difference with the convex non-parametric technology is
that $t\in\{0,1\}^{\ell}$. The FDH technology is non-convex, in
general, but it only postulates the free disposal assumption. Shephard's distance function can also be computed over the FDH
production set by enumeration, see \cite{tve95}. One can also consider mixed approaches combining
both DEA and FDH approaches, see \cite{p05}. The next section presents the parametric viewpoint to estimate the
production set.

\section{Generalized Convexities} \label{Power}

This section presents some concepts of generalized convexity based on power functions providing $\mathbb B$-convex sets with a semi-lattice structure (\cite{bh04}). Then, the algebraic structure underlying the Kolm-Pollack dequantization principle is presented, and finally, their tropical limit sets.

\subsection{Generalized Means and Isomorphic Algebraic Structure}

\noindent This subsection recalls a special class of isomorphism based upon the power functions introduced by \cite{bh04}. They involve a suitable notion of generalized means that was considered by \cite{a70} to analyze social welfare functions. This type of algebraic structure was also studied in \cite{b77} in a convex analysis context. For all $\alpha \in (0,+\infty)$, let $\varphi_\alpha:\Real\longrightarrow
\Real$ be the map defined by:
\begin{equation}\notag
\varphi_{\alpha}(\lambda) =\left\{\begin{array}{ll}
\quad\lambda^{{\alpha}}&\hbox{ if } \lambda \geq 0\\
-|\lambda|^{{\alpha}}&\hbox{ if } \lambda\leq 0
\end{array}
\right.
\end{equation}
For all $\alpha\not=0$, the reciprocal map is defined as
$\varphi_{{\alpha}}^{-1}=\varphi_{\frac{1}{\alpha}}$.
Clearly:
 $(i)$ $\varphi_{\alpha}$ is defined over $\mathbb{R}$;
 $(ii)$ $\varphi_{\alpha}$ is continuous over $\mathbb{R}$;
 $(iii)$ $\varphi_{\alpha}$ is bijective. For all $\lambda, \mu \in \Real$, let us define the operations,
 \begin{align}
\lambda\stackrel{\alpha}{+}\mu &=\varphi_{\alpha}^{-1}\big(\varphi_{\alpha}^{}(\lambda)
+\varphi_{\alpha}^{}(\mu)\big)\notag \\
\lambda\stackrel{\alpha}{\cdot}\mu&=\varphi_{\alpha}^{-1}\left(\varphi_{\alpha}^{}(\lambda)\varphi_{\alpha}^{}(\mu)\right)=\lambda\mu \notag
\end{align}
Then, $(\Real,\stackrel{\alpha}{+}, \cdot) $ is a scalar field, see \cite{b77}. For all vectors
$x=(x_{1},\ldots,x_{d})\in \Real^{d}$, the elementwise power function is defined as:
\begin{equation}\notag
\phi_{\alpha}^{}(x)=\big( {\varphi_\alpha (x_{1}}), \ldots ,
\varphi_\alpha({x_{d}})\big)
\end{equation} 
For $x\in \Real_+^d$, then
$\phi_{\alpha}(x)=\big( {x_{1}}^{\alpha}, \ldots ,
{x_{d}}^{\alpha}\big)=x^{\alpha}.$ The scalar field $(\Real,\stackrel{\alpha}{+}, \cdot) $ is extended to a vector space defining the addition of vectors as follows:
\begin{align}\notag
x\stackrel{\alpha}{+}y=\phi_{\alpha}^{-1}\big(\phi_{\alpha}^{}(x)
+\phi_{\alpha}^{}(y)\big)=(x_1\stackrel{\alpha}{+}y_1,\ldots, x_n\stackrel{\alpha}{+}y_n)
\end{align}
Let $(\Real^n,\stackrel{\alpha}{+}, \cdot) $ denote this vector space. In the following, a subset $C$ of $ \Real^d $ is said to be $\phi_\alpha$-convex if for all $x,y\in C$, all $s,t\in  [0,1]$ and all $\alpha$, $s\stackrel{\alpha}{+}t=1$, implies that $sx\stackrel{\alpha}{+}ty\in C$.  {It follows that a subset $C$ of $\Real^d$ is $\phi_{\alpha}$-convex if, and only if, $\phi_\alpha(C)$ is convex}. Let us consider $A=\{x_{1},\ldots,x_{\ell}\} \subset \Real^{d}$. The $\phi_{\alpha}$-convex hull of the
set $A$ is:
\begin{align}\notag
Co^{\phi_{\alpha}}(A)=\Big\{\sum_{k\in [\ell]}^{\varphi_{\alpha}}t_{k}
 {\cdot}x_{k}:
\sum_{k\in [\ell]}^{\varphi_{\alpha}}t_{k}=1, \ t{\geq} 0\Big\}
\end{align}
If $A\subset \Real_+^d$, then:
\begin{align}\notag
Co^{\phi_{\alpha}}(A)=\Big\{\big(\sum_{k\in
[\ell]}{t_{k}}^{\alpha}{x_{k}\,}^{\alpha}\big)^{\frac{1}{\alpha}} :
\big(\sum_{k\in [\ell]}{t_{k}}^{\alpha}\big)^{\frac{1}{\alpha}}=1, \ t{\geq} 0\Big\}
\end{align}
For the sake of simplicity let us denote $Co^\alpha(A)=Co^{\phi_\alpha}(A)$ for all finite subsets  $A$ of $\Real^d$.

Let us focus now on the case where $\alpha \in (-\infty,0)$. The map
$x\longrightarrow x^{{\alpha}}$ is not defined at point $x=0$. Therefore, in such a case, $\varphi_\alpha$ is not a bijective endomorphism defined on
 $\Real$. Let us denote $K=\{\infty\} \cup \Real \setminus \{0\}$. For all  $\alpha\in (-\infty,0)$ let us consider the map
${  \varphi}_{\alpha}: K\longrightarrow \Real$  defined by:

\begin{equation}\notag
{  \varphi}_{\alpha}(\lambda) =\left\{\begin{array}{ll}
\quad \lambda^{{\alpha}}&\hbox{ if } \lambda > 0\\
-|\lambda|^{{\alpha}}&\hbox{ if } \lambda< 0\\
 \quad 0&\hbox{ if } \lambda = +\infty   \end{array}
                \right.
\end{equation}

\noindent The map ${\varphi}_{\alpha}$
 is then bijective from
$K$ to $\Real$. Thereby, an addition operator $\stackrel{\alpha}{+}$ can be constructed as well as a scalar multiplication $\stackrel{\alpha}{\cdot}$ that are respectively defined for all $\lambda,\mu\in \Real$ as $\lambda \stackrel{\alpha}{+}\mu=\varphi_\alpha^{-1}\big(\varphi_\alpha(\lambda)+\varphi_\alpha^{}(\mu)\big)$ and $\lambda \stackrel{\alpha}{\cdot}\mu=\varphi_\alpha^{-1}\big(\varphi_\alpha(\lambda)\cdot\varphi_\alpha^{}(\mu)\big)$.  Therefore $(K,\stackrel{\alpha}{+},\stackrel{\alpha}{\cdot} )$ has a scalar field structure. In addition, let us introduce the isomorphism $  \phi_\alpha: K^d\longrightarrow \Real^d$, defined by
  $   \phi_\alpha(x_1,\ldots,x_d)=(  \varphi_\alpha (x_1),\ldots,
\varphi_\alpha (x_d))$. For all  $\alpha<0$, let us consider the operations
$\stackrel{{\alpha}}{+}$ and $\stackrel{{\alpha}}{\cdot}$ defined over $K^d$ by:
\begin{align}\notag
x\stackrel{ {\alpha}}{+}y&={  \phi_{\alpha}}^{-1}\left({ 
\phi_{\alpha}}^{}(x)
+{  \phi_{\alpha}}(y)\right)\\
\lambda\stackrel{ {\alpha}}{\cdot }x&={  \phi_{\alpha}}^{-1} \left({
\varphi_{\alpha}}(\lambda). {   \phi_{\alpha}}(x)\right)\notag
\end{align}
Then $\big(K^d, \stackrel{{\alpha}}{+},\stackrel{{\alpha}}{\cdot}\big)$
is a vector space. Notice that if $\alpha<0$  then
$\phi_\alpha =\phi_{-1}(\phi_{|\alpha|})$. 
By extension all these operations are extended to the case $\alpha=0$ using the logarithm function. Hence we defined over $\Real_{++}$ the operation:
\begin{equation}\notag
\lambda\stackrel{0}{+}\mu=\exp (\ln (\lambda)+\ln (\mu)),
\end{equation}
setting $\varphi_0(\lambda)=\ln (\lambda)$ for all positive real numbers $\lambda$. Moreover, let us define $\phi_0(x)=(\varphi_0(x_1),...,\varphi_0(x_n))$ for all $x\in \Real_{++}^d$. By definition, it follows that ${\phi_0}^{-1}(\Real^d)=\Real_{++}^d$.

\begin{center}
{\scriptsize
\unitlength 0.4mm 
\linethickness{0.4pt}
\ifx\plotpoint\undefined\newsavebox{\plotpoint}\fi 
\begin{picture}(219.5,111.25)(0,0)
\put(3,11){\vector(1,0){95.5}}
\put(119.75,11){\vector(1,0){95.5}}
\put(2.75,11){\vector(0,1){92.25}}
\put(119.5,11){\vector(0,1){92.25}}
\put(2.75,103.25){\vector(0,1){0}}
\put(119.5,103.25){\vector(0,1){0}}
\put(0,7.5){\makebox(0,0)[cc]{$0$}}
\put(116.75,7.5){\makebox(0,0)[cc]{$0$}}
\put(102.75,10.75){\makebox(0,0)[cc]{$x$}}
\put(219.5,10.75){\makebox(0,0)[cc]{$x$}}
\put(2.25,111.25){\makebox(0,0)[cc]{$x'$}}
\put(119,111.25){\makebox(0,0)[cc]{$x'$}}
\put(12.575,89.2){\circle*{1.453}}
\put(129.325,89.2){\circle*{1.453}}
\put(92.85,23.55){\circle*{1.3}}
\put(209.6,23.55){\circle*{1.3}}
\put(36.3,46.625){\circle*{1.3}}
\put(193.8,51.625){\circle*{1.3}}
\put(158.55,74.625){\circle*{1.3}}
\qbezier(12.575,89.525)(35.325,88.713)(35.975,46.95)
\qbezier(35.975,46.95)(91.712,46.95)(92.2,23.55)
\qbezier(92.2,23.55)(91.712,89.688)(12.575,90.175)
\put(12.45,80.65){\makebox(0,0)[cc]{$x_1$}}
\put(133.45,93.9){\makebox(0,0)[cc]{$x_1$}}
\put(27.525,41.75){\makebox(0,0)[cc]{$x_2$}}
\put(189.275,60){\makebox(0,0)[cc]{$x_3$}}
\put(154.025,83){\makebox(0,0)[cc]{$x_2$}}
\put(102.6,21.6){\makebox(0,0)[cc]{$x_3$}}
\put(219.35,21.6){\makebox(0,0)[cc]{$x_4$}}
\put(79.75,83.75){\makebox(0,0)[cc]{$Co^\alpha(A)$}}
\put(196.5,83.75){\makebox(0,0)[cc]{$Co^\alpha(A)$}}
\qbezier(129.25,89)(130.75,24.5)(209.25,24)
\qbezier(209.5,23.75)(194,23.5)(193.5,51.25)
\qbezier(193.5,51.25)(160,51.25)(158.5,74.25)
\qbezier(129,89.5)(129.5,74)(157,74.5)
\put(49,-5){\makebox(0,0)[cc]
{ {\bf Figure \ref{Power}.1} Convex hull for $\alpha>1$}}
\put(176.5,-5){\makebox(0,0)[cc]
{ {\bf Figure \ref{Power}.2} Convex hull for $\alpha<1$}}
\end{picture}
}
\end{center}
\bigskip
The limits of these convex hulls may be defined according to semilattice operations $\vee$ and $\wedge$. These operations are intimately linked to the notion of generalized mean. Let us first consider the binary operation $\vee$ defined $\Real^d\times \Real^d$ as,
\begin{equation}\notag
x\vee y=(\max\{x_1,y_1\},\ldots,(\max\{x_d,y_d\})
\end{equation}
Moreover, if $x_1,x_2,\ldots,x_\ell$ are $\ell$ points of $\Real^d$, then:
\begin{equation}\notag
\bigvee_{i\in [\ell]}x_i=\big(\max\{x_{1,1},\ldots,x_{1,\ell} \},\ldots,\max\{x_{1,d},\ldots,x_{\ell,d} \}\big)
\end{equation}
A subset $L$ of $\Real^d$ is an  upper-semilattice if for all $x,y\in L$,
$x\vee y\in L$. In the case of lower-semilattices, for all $x,y\in \Real^d$,
\begin{equation}\notag
x\wedge y=(\min\{x_1,y_1\},\ldots,(\min\{x_d,y_d\})
\end{equation}
Moreover, if $x_1,x_2,\ldots,x_\ell$ are $\ell$ points of $\Real^d$, then:
\begin{equation}\notag
\bigwedge_{i\in [\ell]}^lx_i=\big(\min\{x_{1,1},\ldots,x_{l,\ell} \},\ldots,\min\{x_{1,d},\ldots,z_{\ell,d} \}\big)
\end{equation}
Accordingly, from \cite{bh04}, the $\mathbb B$-convex hull of the finite set $A=\left\{ z_{1},\ldots, z_{\ell}\right\}\subset
\Real_{+}^{d}$ is defined as:
\begin{equation}\notag
\mathbb{B}(A)=
\Big\{\bigvee_{k\in [\ell]} t_{k} z_{k},  {t} \geq 0,
\max_{k\in [\ell]} t_{k}=1\Big\}
\end{equation}
This set is endowed with an upper semi-lattice structure. For
all $A=\left\{ z_{1},\ldots, z_{\ell}\right\}\subset
({\Real}_{++}\cup \{+\infty\})^{d}$, the set,
$$
\mathbb{B}^{-1}(A)=\Big\{\bigwedge_{k\in [\ell]} s_{k} z_{k},
\min_{k \in [\ell]}s_{k}=1,\ \mathbf{s} \geq 0 \Big\}
$$
is called the inverse
$\mathbb{B}$-convex hull of $A$. These sets were considered in \cite{adilYe,bl11} who established some related properties. These sets can be viewed as the limit of the generalized $\phi_\alpha$-convex hull and are based on a lower-semilattice structure.  

The  lower [upper] Painlev\'e-Kuratowski  limit of a sequence of sets
$\{E_n\}_{n\in\mathbb{N}}$ is denoted $Li_{n\to\infty}E_n$ [$Ls_{n\to\infty}A_n$]. The lower Painlevé-Kuratowski is the set of all the points
 $p$ such that there exists a sequence $\{p_n\}_{n\in \mathbb N}$ such that $p_n\in E_n$ for any $n$ and $\lim_{n\longrightarrow \infty}{p_n}=p$. The upper limit is the set of all the points $p$ such that here exists a subsequence $\{p_{n_k}\}_{k\in \mathbb N}$ such that $p_{n_k}\in E_{n_k}$ for any $k$ and $\lim_{k\longrightarrow \infty}{p_{n_k}}=p$.  By definition $\Linf_{n\to\infty}E_n \subset \Ls_{n\to\infty}E_n$. If the upper and lower limits are identical the sequence of sets $\{E_n\}_{n\in\mathbb{N}}$ has a Painlevé-Kuratowski
 limit and then,
$\Ls_{n\to\infty}E_n =E = \Linf_{n\to\infty}E_n=\Lim_{n\to\infty}E_n$.  

Let $\{\alpha_k\}_{k\in \mathbb N}$ be a sequence of real numbers such that $\lim_{k \longrightarrow\infty}\alpha_k=+\infty$, from \cite{bh04}, we have:
\begin{equation}\notag
\mathbb B(A)=\Lim_{k\longrightarrow \infty}Co^{\alpha_k}(A):=Co^\infty(A)
\end{equation}
Moreover, from \cite{adilYe}, if $A\subset \Real_{++}^d$, and if $ \lim_{k \longrightarrow \infty}\alpha_k=-\infty$, then:
\begin{equation}\notag
\mathbb B^{-1}(A)=\Lim_{k\longrightarrow \infty}Co^{\alpha_k}(A):=Co^{-\infty}(A) 
\end{equation}
Those $\mathbb B$-convex structures have been shown to have relevant properties for production technologies, see \cite{abm17,abs17}.

\subsection{Kolm-Pollack Form: Dequantization}

In this subsection a special generalized mean based on the composition of logarithm and exponential functions is investigated, see \cite{HLP}, in order to provide new algebraic structures for production technologies. 

The so-called Kolm-Pollack form is intimately linked to the generalized mean and to some concepts of tropical geometry, and in particular, to the dequantization principle analyzed in \cite{ms}. Let us first define for all $\alpha\in \Real\backslash \{0\}$, the binary operation $\stackrel{\alpha}{\oplus}: \Real\times \Real\longrightarrow \Real$ as,
\begin{equation}\notag
\lambda \stackrel{\alpha}{\oplus}\mu=\frac{1}{\alpha}\ln \big(e^{\alpha \lambda} + e^{\alpha \mu}\big)
\end{equation}
This binary operation has the form $\lambda \stackrel{\alpha}{\oplus}\mu=\varsigma_{\alpha}^{-1}\big (\varsigma_\alpha (\lambda)+\varsigma_\alpha(\mu)\big)$ where  {$\varsigma_\alpha=\varphi_\alpha\circ \exp$.}   The product related to this operation is the standard addition:
\begin{equation}\notag
\lambda \stackrel{\alpha}{\otimes}\mu=\frac{1}{\alpha}\ln \big(e^{\alpha \lambda} \cdot e^{\alpha \mu}\big)=\lambda +\mu
\end{equation}
For all $\alpha \not=0$, let   $\mathbb M_\alpha=\Real\cup\{-\infty\}$ if $\alpha >0$ and $\mathbb M_\alpha=\Real\cup\{+\infty\}$ if $\alpha <0$. By construction $(\mathbb M_\alpha,\stackrel{\alpha}{\oplus}, +)$ is an idempotent semi-ring. These operations can be extended to $ \mathbb M_\alpha ^{d}$: 
\begin{equation}\notag
x\stackrel{\alpha}{\oplus}y=(x_1\stackrel{\alpha}{\oplus }y_1,\ldots, x_d\stackrel{\alpha}{\oplus }y_d)
\end{equation}
For all $x\in 
\mathbb M_\alpha^d$, the generalized mean  is given by, 
\begin{equation}\notag
{\stackrel{\alpha}{\bigoplus_{i\in [d]}}}\  x_i=\frac{1}{\alpha}\ln \Big(\sum_{i\in [d]}e^{\alpha x_i}\Big)
\end{equation}
Let $\sigma_\alpha: \mathbb M_\alpha ^d\longrightarrow \mathbb M_\alpha^d$ be the map defined for all $(x_1,\ldots,x_d)\in \Real^d$ as $\sigma_\alpha(x_1,\ldots,x_d)=(\varsigma_\alpha(x_1),\ldots,\varsigma_\alpha(x_d))$. Let $\big(\mathbb M_\alpha^{d}, \stackrel{\alpha}{\oplus},+\big)$  denote this semi-ring. 

\begin{defn}
 {For all $\alpha \neq 0$}, a subset $M$ of $\mathbb M_\alpha^d$ is $\psi_\alpha$-convex if for all $x,y\in M$ and all $s,t\in  -\Real$, $s\stackrel{\alpha}{\oplus}t=0$ implies that $(s1\!\!1_d \otimes x)\stackrel{\alpha}{\oplus}(t1\!\!1_d \otimes y)\in M$. It follows that a subset $M$ of $\mathbb M_\alpha^d$ is $\psi_\alpha$-convex  {if, and only if, $\sigma_\alpha(M)$ is convex}.
\end{defn}
Let us consider $A=\{x_{1},\ldots,x_{\ell}\} \subset \mathbb M^{d}_\alpha$, then the $\sigma_{\alpha}$-convex hull of the
set $A$ is given by,
\begin{align}\notag
Co^{\sigma_{\alpha}}(A)=\Big\{ \stackrel{\alpha}{\bigoplus _{k\in
[\ell]}}t_{k}\otimes x_{k} : \ 
\stackrel{\alpha}{\bigoplus _{k\in
[\ell]}}t_{k} =0,  t\in \mathbb M_\alpha^\ell \Big\}
\end{align}
Equivalently:
\begin{align}\notag
Co^{\sigma_{\alpha}}(A)=\Big\{\frac{1}{\alpha}\mathbf{ln}\big(\stackrel{\alpha}{\sum _{k\in
[\ell]}}\mathbf{e}^{\alpha (t_{k}1\!\!1_d+x_k)}\big) :
\frac{1}{\alpha}\ln\big(\sum_{k\in [\ell]}e^{\alpha t_{k}} \big)=0, \ t\in \mathbb M_\alpha^\ell \Big\}
\end{align}

\subsection{Tropical Limit Sets}\label{LimSet}

By moving the value of $\alpha$, the shape of the generalized convex hull is modified, and this allows new properties to be designed for production technologies. Moreover, taking limit values for $\alpha$ provides particular convex hulls related to tropical geometries. The limit of the simplex is analyzed when   {$\alpha\longrightarrow \alpha_{0}\in \Real \cup
\{-\infty, +\infty\} \setminus \{0\}$.} This implies a geometric deformation of the $\sigma_{\alpha}$-convex hull.

For all finite and nonempty set $A$ contained in $\Real^d$,
$Co^{\sigma_\alpha}(A)$ belongs to ${\cal K}(\Real^d)$, the space of nonempty
compact subsets of $\Real^d$, which is metrizable by the
Hausdorff metric:
$$
D_{H}(C_1, C_2) = \inf\Big\{\varepsilon > 0: C_1\subset\bigcup_{x\in C_{2}}B(x,
\varepsilon), \hbox{ and } C_2\subset\bigcup_{x\in C_{1}}B(x,
\varepsilon)\Big\}
$$
The so-called dequantization principle   consists in considering the limit of the algebraic structure $(\mathbb M_\alpha,\stackrel{\alpha}{\oplus}, +)$ when
$\alpha\longrightarrow \infty$. By convention let us define $\mathbb M_\infty=\Real\cup\{-\infty\}$ and $ \mathbb M_{-\infty}=\Real\cup\{+\infty\}$. For all $\lambda,\mu\in \mathbb M_{\infty}$, we have:
\begin{equation}\notag
\lim_{\alpha\longrightarrow \infty}\lambda \stackrel{\alpha}{\oplus}\mu=\max\{\lambda,\mu\}:=\lambda \stackrel{\infty}{\oplus}\mu
\end{equation}
 {For all $\lambda,\mu\in \mathbb M_{-\infty}$, we have:}
\begin{equation}\notag
\lim_{\alpha\longrightarrow -\infty}\lambda \stackrel{\alpha}{\oplus}\mu=\min\{\lambda,\mu\}:=\lambda \stackrel{-\infty}{\oplus}\mu
\end{equation}
It follows that $(\mathbb M_{\infty},\stackrel{\infty}{\oplus}, + )$ and $(\mathbb M_{-\infty},\stackrel{-\infty}{\oplus}, + )$ respectively define an idempotent semi-ring over 
$\mathbb M_{\infty}$ and $\mathbb M_{-\infty}$.  Let us consider $A=\{x_{1},\ldots,x_{\ell}\} \subset \mathbb M_\infty^{d}$ we have:
\begin{align}\notag
\lim_{\alpha\longrightarrow \infty}  \stackrel{\alpha}{\bigoplus _{k\in
[\ell]}}t_{k}\otimes x_{k} =\bigvee_{k\in [\ell]}  (t_{k}1\!\!1_d+x_k) =\stackrel{\infty}{\bigoplus _{k\in
[\ell]}}t_{k}\otimes x_{k}
\end{align}
If $A=\{x_{1},\ldots,x_{\ell}\} \subset \mathbb M_{-\infty}^{d}$ we have:
\begin{align}\notag
\lim_{\alpha\longrightarrow - \infty}  \stackrel{\alpha}{\bigoplus_{k\in
[\ell]}}t_{k}\otimes x_{k} =\bigwedge_{k\in [\ell]}  (t_{k}1\!\!1_d+x_k) =\stackrel{-\infty}{\bigoplus _{k\in
[\ell]}}t_{k}\otimes x_{k}
\end{align}

\begin{defn}
Given a subset $M$ of $\mathbb M_{\infty}^d$, $M$ is said to be Max-Plus convex if, for all $s,t\leq 0$ and all $x,y\in M$, then $(s\otimes x)\stackrel{ \infty}{\oplus}(t\otimes y)\in M$. 
\end{defn}
Let us consider $A=\{x_{1},\ldots,x_{\ell}\} \subset \mathbb M_\infty^{d}$. The Max-Plus convex hull of the
set $A$ is:
\begin{align}\notag
Co^{\sigma_{\infty}}(A)&=\Big\{ \stackrel{\infty}{\bigoplus _{k\in
[\ell]}}t_{k}\otimes x_{k} :
\ \stackrel{\infty}{\bigoplus _{k\in
[\ell]}}t_{k} =0,  t\in \mathbb M_\infty^\ell \Big\}\\
&=\Big\{  {\bigvee _{k\in
[\ell]}}  t_{k}1\!\!1_d+x_k  :
\max_{k\in [\ell]}t_k=0, \ t\in \mathbb M_{\infty}^\ell \Big\} \notag
\end{align}
 Let us consider $A=\{x_{1},\ldots,x_{\ell}\} \subset \mathbb M_{-\infty}^{d}$. The Min-Plus convex hull of the
set $A$ is:
\begin{align}\notag
Co^{\sigma_{-\infty}}(A)&=\Big\{ \stackrel{-\infty}{\bigoplus _{k\in
[\ell]}}t_{k}\otimes x_{k} : \
\stackrel{-\infty}{\bigoplus _{k\in
[\ell]}}t_{k} =0,  t\in \mathbb M_{-\infty}^\ell \Big\}\\
&=\Big\{  {\bigwedge _{k\in
[\ell]}}  t_{k}1\!\!1_d+x_k  :
\min_{k\in [\ell]}t_k=0, \ t\in \mathbb M_{-\infty}^\ell \Big\} \notag 
\end{align}
 
By construction interesting properties arise between the geometry based on the power functions and that based on exponential functions. 

\begin{lem}
A subset $M$ of $\mathbb M_{\infty}^d $ is Max-Plus convex if, and only if, $\mathbf{exp}(M)$ is $\mathbb B$-convex. Symmetrically 
a subset $M$ of $\mathbb M_{-\infty}^d $ is Min-Plus convex if, and only if, $ \mathbf{exp}\circ \phi_{-1}(M)$ is inverse $\mathbb B$-convex.
\end{lem}

Let us now investigate the Painlevé-Kuratowski limits of the Max-Plus and the Min-Plus convex hulls.

\begin{prop}\label{LimTRopConv}
 {Suppose that $A\subset  \mathbb M^d_{\infty}$}. Let $\{\alpha_k\}_{k\in \mathbb N}$ be a sequence of real numbers such that $\lim_{k\longrightarrow \infty}\alpha_k=+\infty. $ Then:
 
\begin{equation*}
\Lim_{k\longrightarrow \infty}Co^{\sigma_{\alpha_k}}(A)=Co^{\sigma_{\infty}}(A)
\end{equation*}
Suppose now that  {$A\subset  \mathbb M^d_{-\infty}$} and $\lim_{k\longrightarrow \infty}\alpha_k=-\infty$. Then:
\begin{equation*}
\Lim_{k\longrightarrow -\infty}Co^{\sigma_{\alpha_k}}(A)=Co^{\sigma_{-\infty}}(A)
\end{equation*}
\end{prop}
{\bf Proof:} We have:
\begin{align}\notag
Co^{\sigma_{\alpha_k}}(A)&= \frac{1}{\alpha_k}\mathbf{ln}\Big(Co(\mathbf e^{\alpha_k A}  )\Big)=\mathbf{ln}\Big(Co^{\alpha_k}(\mathbf e^{ A})\Big)\end{align}
 However $\mathbf e^{ A}\in \Real_+^d.$ From \cite{bh04}, $\Lim_{k\longrightarrow \infty}Co^{\alpha_k}(\mathbf e^{ A})=Co^{\infty}(\mathbf e^{ A})$. By definition,
 $$
 Co^{\infty}(\mathbf e^{ A})=\Big\{\bigvee_{k\in [\ell]}t_k\mathbf e^{z_k}: \max_{k\in [\ell]} t_k=1, t\geq 0\Big\}\subset \Real_{++}^d
 $$ 
 Setting $s_k=\ln t_k$ for each $k$, with the convention $-\infty=\ln (0)$, then: 
$$
Co^{\infty}(\mathbf e^{ A})=\Big\{\bigvee_{k\in [\ell]} \mathbf e^{s_k \un_d+ z_k}: \max_{k\in [\ell]} s_k=0, s\in \mathbb M_{\infty}^\ell\Big\}
$$
The map $\mathbf{ln}: \Real_{++}^d\mapsto \Real^d$ is continuous over  $\Real_{++}^d$, then it follows that:
$$\Lim_{k\longrightarrow \infty}\mathbf{ln}\Big(Co^{\alpha_k}(\mathbf e^{ A})\Big)=\mathbf{ln}\Big(\Lim_{k\longrightarrow \infty}Co^{\alpha_k}(\mathbf e^{ A})\Big)=\mathbf{ln}\Big(Co^{\infty}(\mathbf e^{ A})\Big)
$$
However:
\begin{align*}
\mathbf{ln}\Big(Co^{\infty}(\mathbf e^{ A})\Big)&=\mathbf{ln}\Big(\Big\{\bigvee_{k\in [\ell]} \mathbf e^{s_k \un_d+ z_k}: \max_{k\in [\ell]} s_k=0, s\in \mathbb M_{\infty}^\ell\Big\}\Big)\\&=\Big\{  {\bigvee _{k\in
[\ell]}}  s_{k}1\!\!1_d+x_k  :
\max_{k\in [\ell]}s_k=0, s\in \mathbb M_{\infty}^\ell \Big\}\\&=Co^{\sigma_{\infty}}(A)
\end{align*}
We deduce the first part of the statement. The proof of the second part is similar. $\Box$\\

\begin{center}{\scriptsize
\unitlength 0.4mm 
\linethickness{0.4pt}
\ifx\plotpoint\undefined\newsavebox{\plotpoint}\fi 
 
\unitlength 0.4mm 
\linethickness{0.4pt}
\ifx\plotpoint\undefined\newsavebox{\plotpoint}\fi 

\begin{picture}(229.25,111.75)(0,0)
\put(56.25,91){\circle*{1.803}}
\put(168.5,78.5){\circle*{1.803}}
\put(131.5,59.75){\circle*{1.803}}
\put(51.25,103.75){\circle*{1.803}}
\put(197,103.75){\circle*{1.803}}
\put(54.5,35.75){\circle*{1.803}}
\put(78.25,80.25){\circle*{1.803}}
\put(224,80.25){\circle*{1.803}}
\put(75.75,67.5){\circle*{1.803}}
\put(202.5,45.25){\circle*{1.803}}
\put(45.25,27){\circle*{1.803}}
\put(14.75,27){\circle*{1.803}}
\put(56.25,66.25){\circle*{1.803}}
\put(80.5,36.75){\circle*{1.803}}
\put(9,106.25){\circle*{1.803}}
\put(9,77.75){\circle*{1.803}}
\put(2.25,50.25){\circle*{1.581}}
\put(188.75,23.25){\circle*{1.581}}
\put(13.5,51){\line(-1,0){11}}
\put(51.25,104){\line(1,0){26.5}}
\put(197,104){\line(1,0){26.5}}
\put(54.5,36){\line(1,0){26.5}}
\put(77.75,104){\line(0,-1){22.5}}
\put(223.5,104){\line(0,-1){22.5}}
\put(45,39.75){\line(0,-1){13.5}}
\multiput(15.25,27.75)(.03500859107,.03371993127){1164}{\line(1,0){.03500859107}}
\put(9,106.5){\line(0,-1){27.5}}
\put(56.25,95.25){\makebox(0,0)[cc]{$x_1$}}
\put(168.5,82.75){\makebox(0,0)[cc]{$x_3$}}
\put(131.5,64){\makebox(0,0)[cc]{$x_1$}}
\put(8.5,111.75){\makebox(0,0)[cc]{$x_1$}}
\put(50.5,109){\makebox(0,0)[cc]{$x_1$}}
\put(196.25,109){\makebox(0,0)[cc]{$x_1$}}
\put(54.5,72.25){\makebox(0,0)[cc]{$x_1$}}
\put(44.75,20.5){\makebox(0,0)[cc]{$x_1$}}
\put(53.5,29.5){\makebox(0,0)[cc]{$x_1$}}
\put(0,45.25){\makebox(0,0)[cc]{$x_2$}}
\put(186.5,18.25){\makebox(0,0)[cc]{$x_5$}}
\put(9.5,72.75){\makebox(0,0)[cc]{$x_2$}}
\put(83.5,80){\makebox(0,0)[cc]{$x_2$}}
\put(229.25,80){\makebox(0,0)[cc]{$x_2$}}
\put(81.75,68.5){\makebox(0,0)[cc]{$x_2$}}
\put(208.5,46.25){\makebox(0,0)[cc]{$x_4$}}
\put(86.75,34.75){\makebox(0,0)[cc]{$x_2$}}
\put(10.25,20.5){\makebox(0,0)[cc]{$x_2$}}
\multiput(223.5,80.5)(-.0337301587,-.0341269841){630}{\line(0,-1){.0341269841}}
\put(202.25,59){\line(0,-1){13.75}}
\put(188.5,60){\line(0,-1){36.25}}
\put(57,0){\makebox(0,0)[cc]
{{{\bf Figure \ref{Power}.3:}  Max-Plus convex hull of two points.\ \ }}}
\put(191.5,0){\makebox(0,0)[cc]{{{\bf Figure \ref{Power}.4:} Max-Plus-polytope.}}}
\multiput(13.75,51.25)(.03586678053,.03373185312){1171}{\line(1,0){.03586678053}}
\multiput(45.25,39.75)(.03675577157,.0337181045){823}{\line(1,0){.03675577157}}
\multiput(196.75,103.5)(-.03373579545,-.03586647727){704}{\line(0,-1){.03586647727}}
\put(173,78.25){\line(-1,0){5}}
\put(188.5,60){\line(-1,0){44.5}}
\put(188.45,60){\line(-1,0){57.85}}
\multiput(168,78.5)(-.0337338262,-.0341959335){541}{\line(0,-1){.0341959335}}
\multiput(201.5,45.5)(-.0337150127,-.0388040712){393}{\line(0,-1){.0388040712}}
\end{picture}

}
\end{center}

\begin{center}
{\scriptsize 
\unitlength 0.4mm 
\linethickness{0.4pt}
\ifx\plotpoint\undefined\newsavebox{\plotpoint}\fi 
\begin{picture}(233.901,114.5)(0,0)
\put(30.5,44){\circle*{1.803}}
\put(196,53){\circle*{1.803}}
\put(233,71.75){\circle*{1.803}}
\put(35.5,31.25){\circle*{1.803}}
\put(167.5,27.75){\circle*{1.803}}
\put(32.25,99.25){\circle*{1.803}}
\put(8.5,54.75){\circle*{1.803}}
\put(140.5,51.25){\circle*{1.803}}
\put(11,67.5){\circle*{1.803}}
\put(162,86.25){\circle*{1.803}}
\put(41.5,108){\circle*{1.803}}
\put(72,108){\circle*{1.803}}
\put(30.5,68.75){\circle*{1.803}}
\put(6.25,98.25){\circle*{1.803}}
\put(77.75,28.75){\circle*{1.803}}
\put(77.75,57.25){\circle*{1.803}}
\put(84.5,84.75){\circle*{1.581}}
\put(175.75,108.25){\circle*{1.581}}
\put(73.25,84){\line(1,0){11}}
\put(35.5,31){\line(-1,0){26.5}}
\put(167.5,27.5){\line(-1,0){26.5}}
\put(32.25,99){\line(-1,0){26.5}}
\put(9,31){\line(0,1){22.5}}
\put(141,27.5){\line(0,1){22.5}}
\put(41.75,95.25){\line(0,1){13.5}}
\multiput(71.5,107.25)(-.03500859107,-.03371993127){1164}{\line(-1,0){.03500859107}}
\put(77.75,28.5){\line(0,1){27.5}}
\put(30.5,39.75){\makebox(0,0)[]{$x_1$}}
\put(196,48.75){\makebox(0,0)[]{$x_3$}}
\put(233,67.5){\makebox(0,0)[]{$x_1$}}
\put(78.25,23.25){\makebox(0,0)[]{$x_1$}}
\put(36.25,26){\makebox(0,0)[]{$x_1$}}
\put(168.25,22.5){\makebox(0,0)[]{$x_1$}}
\put(32.25,62.75){\makebox(0,0)[]{$x_1$}}
\put(42,114.5){\makebox(0,0)[]{$x_1$}}
\put(33.25,105.5){\makebox(0,0)[]{$x_1$}}
\put(86.75,89.75){\makebox(0,0)[]{$x_2$}}
\put(178,113.25){\makebox(0,0)[]{$x_5$}}
\put(77.25,62.25){\makebox(0,0)[]{$x_2$}}
\put(3.25,55){\makebox(0,0)[]{$x_2$}}
\put(135.25,51.5){\makebox(0,0)[]{$x_2$}}
\put(5,66.5){\makebox(0,0)[]{$x_2$}}
\put(156,85.25){\makebox(0,0)[]{$x_4$}}
\put(0,100.25){\makebox(0,0)[]{$x_2$}}
\put(76.5,114.5){\makebox(0,0)[]{$x_2$}}
\multiput(141,51)(.0337301587,.0341269841){630}{\line(0,1){.0341269841}}
\put(162.25,72.5){\line(0,1){13.75}}
\put(176,71.5){\line(0,1){36.25}}
\put(46,0){\makebox(0,0)[cc]
{{{\bf Figure \ref{Power}.5:}  Min-Plus convex hull of two points.}}}
\put(191,0){\makebox(0,0)[cc]{{{\bf Figure \ref{Power}.6:} Min-Plus-polytope.}}}
\multiput(73,83.75)(-.03586678053,-.03373185312){1171}{\line(-1,0){.03586678053}}
\multiput(41.5,95.25)(-.03675577157,-.0337181045){823}{\line(-1,0){.03675577157}}
\multiput(167.75,28)(.03373579545,.03586647727){704}{\line(0,1){.03586647727}}
\put(191.5,53.25){\line(1,0){5}}
\put(176,71.5){\line(1,0){44.5}}
\put(176.05,71.5){\line(1,0){57.85}}
\multiput(196.5,53)(.0337338262,.0341959335){541}{\line(0,1){.0341959335}}
\multiput(163,86)(.0337150127,.0388040712){393}{\line(0,1){.0388040712}}
\end{picture}
}
\end{center}

\begin{center}{\scriptsize 
\noindent {{\bf Table. 1}: Typology of Limit Sets}\\*
\begin{tabular}
{|p{2 cm}||p{10cm} |}
\hline   &  \vspace{0,15cm}
{ Convex Hull and Convex Hull in Limit}\vspace{0,25cm}
 \\
\hline \vspace{0,3cm}$\alpha \in \Real \setminus \{0\}$ \vspace{0,25cm}& \vspace{0,15cm}$Co^{\sigma_{\alpha}}(A)=\Big\{\frac{1}{\alpha}\mathbf{ln}\big(\stackrel{\alpha}{\sum _{k\in
[\ell]}}\mathbf{e}^{\alpha (t_{k}1\!\!1_d+x_k)}\big) :
\frac{1}{\alpha}\ln\big(\sum_{k\in [\ell]}e^{\alpha t_{k}} \big)=0, \ t\in \mathbb M_\alpha^\ell \Big\}$
\vspace{0,25cm}\\
\hline \vspace{0,3cm}$\alpha_k\longrightarrow +\infty$\vspace{0,25cm} &
\vspace{0,15cm}$Co^{\sigma_{\infty}}(A) =\Big\{  {\bigvee _{k\in
[\ell]}}  t_{k}1\!\!1_d+x_k  :
\max_{i\in [\ell]}t_i=0, \ t\in \mathbb M_{ \infty}^\ell \Big\}$
\vspace{0,25cm}
\\
\hline \vspace{0,25cm}$\alpha_k\longrightarrow -\infty$ \vspace{0,25cm}&
\vspace{0,15cm} $ Co^{\sigma_{-\infty}}(A) =\Big\{  {\bigwedge _{k\in
[\ell]}}  t_{k}1\!\!1_d+x_k  :
\min_{i\in [\ell]}t_i=0, \ t\in \mathbb M_{-\infty}^\ell \Big\}$
\vspace{0,25cm}\\
\hline
\end{tabular}}
\end{center}

In what follows, it will be useful to consider the  generalized $\sigma_\alpha$-simplexes defined as: 

\begin{equation}\notag
\Delta_\ell^{(\sigma_\alpha)}=\Big \{ {t\in \big(\sigma_\alpha(\Real_+)\big)^\ell}: 
\frac{1}{\alpha}\ln\big(\sum_{k\in [\ell]}e^{\alpha t_{k}} \big)=0, \ t\in \mathbb M_\alpha^\ell \Big \}
\end{equation}
From the $\sigma_\alpha$-simplex, the following limiting cases are deduced:
$$
\Delta_\ell^{(\sigma_\infty)}=\Big \{t \in  \Real_+^\ell:   \max\limits_{k\in [\ell]} {t_k}=0, t_k\geq 0\Big \}
$$
$$
\Delta_\ell^{(\sigma_{-\infty})}=\Big \{t\in  \Big(\Real_{++}\cup\{+\infty\}\Big)^\ell: \min\limits_{k\in [\ell]} t_k=0, t_k\in [1,\infty)\Big \}
$$

\begin{center}
{\scriptsize 
\unitlength 0.4mm 
\linethickness{0.4pt}
\ifx\plotpoint\undefined\newsavebox{\plotpoint}\fi 
\begin{picture}(218.692,111.534)(0,0)
\put(192.901,51.284){\circle*{1.803}}
\put(37.901,81.034){\circle*{1.803}}
\put(.902,62.284){\circle*{1.803}}
\put(164.401,26.034){\circle*{1.803}}
\put(66.402,106.284){\circle*{1.803}}
\put(137.401,49.534){\circle*{1.803}}
\put(93.401,82.784){\circle*{1.803}}
\put(158.901,84.534){\circle*{1.803}}
\put(71.902,47.784){\circle*{1.803}}
\put(172.651,106.534){\circle*{1.581}}
\put(217.901,69.534){\circle*{1.581}}
\put(58.152,25.784){\circle*{1.581}}
\put(164.401,25.784){\line(-1,0){26.5}}
\put(66.402,106.534){\line(1,0){26.5}}
\put(137.901,25.784){\line(0,1){22.5}}
\put(92.901,106.534){\line(0,-1){22.5}}
\put(192.901,47.034){\makebox(0,0)[]{$x_3$}}
\put(37.901,85.284){\makebox(0,0)[]{$x_3$}}
\put(.902,66.534){\makebox(0,0)[]{$x_1$}}
\put(165.151,20.784){\makebox(0,0)[]{$x_1$}}
\put(65.652,111.534){\makebox(0,0)[]{$x_1$}}
\put(174.901,111.534){\makebox(0,0)[]{$x_5$}}
\put(55.902,20.784){\makebox(0,0)[]{$x_5$}}
\put(132.151,49.784){\makebox(0,0)[]{$x_2$}}
\put(98.651,82.534){\makebox(0,0)[]{$x_2$}}
\put(152.901,83.534){\makebox(0,0)[]{$x_4$}}
\put(77.902,48.784){\makebox(0,0)[]{$x_4$}}
\multiput(137.901,49.284)(.0843253968,.0853174603){252}{\line(0,1){.0853174603}}
\multiput(92.901,83.034)(-.0843253968,-.0853174603){252}{\line(0,-1){.0853174603}}
\put(159.151,70.784){\line(0,1){13.75}}
\put(71.652,61.534){\line(0,-1){13.75}}
\put(172.901,69.784){\line(0,1){36.25}}
\put(57.902,62.534){\line(0,-1){36.25}}
\put(36.328,2.25){\makebox(0,0)[cc]
{{{\bf Figure \ref{Power}.8:} Max-Plus-polytope in Limit.\ \ \ \ \ \ \ \ \ \  }}}
\multiput(164.651,26.284)(.0842198582,.0895390071){282}{\line(0,1){.0895390071}}
\multiput(66.152,106.034)(-.0842198582,-.0895390071){282}{\line(0,-1){.0895390071}}
\put(188.401,51.534){\line(1,0){5}}
\put(42.401,80.784){\line(-1,0){5}}
\put(172.901,69.784){\line(1,0){44.5}}
\put(57.902,62.534){\line(-1,0){44.5}}
\put(57.852,62.534){\line(-1,0){57.85}}
\multiput(193.401,51.284)(.084101382,.085253456){217}{\line(0,1){.085253456}}
\multiput(37.401,81.034)(-.084101382,-.085253456){217}{\line(0,-1){.085253456}}
\multiput(159.901,84.284)(.083860759,.096518987){158}{\line(0,1){.096518987}}
\multiput(70.902,48.034)(-.083860759,-.096518987){158}{\line(0,-1){.096518987}}
\put(176.401,75.284){\makebox(0,0)[]{$x_1$}}
\put(210.901,75.284){\makebox(0,0)[]{$x_2$}}
\put(168.578,2.25){\makebox(0,0)[cc]
{{{\bf Figure \ref{Power}.8:} Min-Plus-polytope in Limit. }}}
\qbezier(38.151,81.034)(19.151,63.034)(1.151,63.034)
\qbezier(164.151,25.784)(188.276,51.534)(192.901,51.284)
\qbezier(66.651,106.534)(42.526,80.784)(37.901,81.034)
\qbezier(164.401,25.534)(137.776,25.784)(137.651,49.034)
\qbezier(66.401,106.784)(93.026,106.534)(93.151,83.284)
\qbezier(137.401,49.784)(159.276,71.034)(158.651,84.284)
\qbezier(93.401,82.534)(71.526,61.284)(72.151,48.034)
\qbezier(172.401,105.534)(172.401,98.284)(158.401,84.034)
\qbezier(58.401,26.784)(58.401,34.034)(72.401,48.284)
\qbezier(.651,62.534)(57.776,62.784)(57.401,26.034)
\qbezier(193.324,51.297)(210.878,68.85)(217.92,69.166)
\qbezier(173.563,106.586)(173.353,69.376)(217.71,70.007)
\end{picture}
}
\end{center}
\bigskip

\section{Quantized Convexity and Non-Parametric Production Technology}\label{Quantized}

Subsection 1 introduces input and output distance functions in $\sigma_\alpha$-convex spaces. Subsection 2 extends these distance functions to non-parametric production models. Subsection 3 is devoted to the dual properties of these distance functions. 

\subsection{Distance Functions}
Paralleling the usual multiplication of vectors by scalar numbers,
 {for all $s\in \mathbb M_{\infty}$ and all $z\in \mathbb M_{\infty}^d$ the}
Max-Plus multiplication by a scalar number is defined by:

\begin{equation}\notag
s\otimes z :=(s \otimes z_1,\ldots,s \otimes z_d)=(s + z_1,\ldots,s +
z_d)=z+s\un_d.
\end{equation}
Accordingly, the input translation function can then  be defined as,

\begin{equation}\notag
\mathbb D_{\mathrm{in}}(x,y)=\sup \{\delta\in \Real: (-\delta )\otimes x\in L(y)\}.
\end{equation}
On the other hand, the output translation function is,

\begin{equation}\notag
\mathbb D_{\mathrm{out}}(x,y)=\sup \{\delta\in \Real: \delta \otimes y\in P(y)\}.
\end{equation}

The translation distance function can be viewed as a restricted case of  the {  topical functions} introduced in \cite{gunkea} (see \cite{rubising} for related topics). A function $f: \Real^m\longrightarrow \Real\cup\{-\infty,+\infty\}$ is called  {  topical} if this function is weakly monotonic with respect to the usual partial order defined over $\Real^m$ and satisfies translation homotheticity ($f(x+\alpha \un)=f(x)+\alpha$ for all $x \in \Real^m$ and all $\alpha \in \Real$). It follows that for all $y\in \Real_+^n$ the map $x\mapsto \mathbb D_{\mathrm{in}}(x,y)$ satisfies the translation  property of topical functions when $x+\alpha \un_m \in \Real_+^m$. The input oriented translation function is also related to the { nonlinear scalarization function}  defined in \cite{Tamm83} and \cite{pasco} (see also \cite{CHY} and Definition 1.40, p.13). More recently, an extensive foundation for the application of the concept of scalarization has been provided in \cite{Tamm20}. Making obvious changes of variables, it is also related to the nonlinear functional introduced in \cite{GRTZ03} (see Definition 2.23, p.39).

Therefore, the Graph-translation homotheticity property can be defined as:

\medskip

\noindent T4: ($\delta \otimes x,\delta \otimes y)\in T$ for all $\delta \in \Real$ such that $(\delta \otimes x,\delta
\otimes y)\geq 0$.

\medskip

In the following, it is shown that the distance functions $\mathbb D_{\mathrm{in}}$ and $\mathbb D_{\mathrm{out}}$ can be computed over transformed technologies $\mathbf{e}^{\alpha T}$ with Farrell input and output efficiency scores ($E_{\mathrm{in}}$ and $E_{\mathrm{out}}$, respectively).  
 
\begin{prop}\label{Distance}
Suppose that $T\in\mathcal T$. Then:
\newline $(i)$ $\mathbb D_{\mathrm{in}}(x,y,T)=-\frac{1}{\alpha}\ln\big(E_{\mathrm{in}}(\mathbf e^{\alpha x}, \mathbf e^{\alpha y }, \mathbf{e}^{\alpha T})\big) $
\newline $(ii)$ $\mathbb D_{\mathrm{out}}(x,y,T)=\frac{1}{\alpha}\ln\big(E_{\mathrm{out}}(\mathbf e^{\alpha x}, \mathbf e^{\alpha y}, \mathbf{e}^{\alpha T})\big)$

 \end{prop}
 {\bf Proof:} $(i)$ We have,
$$
\mathbb D_{\mathrm{in}}(x,y,T)=\sup\{\delta: (x-\delta \un_m, y)\in T)\}   =\sup\{\delta: (e^{-\alpha \delta }\mathbf e^{\alpha x}, \mathbf e^{ \alpha y})\in \mathbf e^{\alpha T)}\} 
$$
Setting $ \lambda =e^{-\alpha \delta}$ we deduce that $\delta =-\frac{1}{\alpha}\ln (\lambda)$. Therefore:
$$   
\mathbb D_{\mathrm{in}}(x,y,T)=\sup \big \{-\frac{1}{\alpha}\ln (\lambda): (\lambda \mathbf e^{\alpha x}, \mathbf e^{ \alpha y})\in \mathbf e^{\alpha T)}\big\}  = -\frac{1}{\alpha}\ln\Big( \inf \big \{\lambda: (\lambda \mathbf e^{\alpha x}, \mathbf e^{ \alpha y})\in \mathbf e^{\alpha T)}\big\}\Big)
$$
which proves the first part of the statement. $(ii)$ We have:
\begin{equation}
\mathbb D_{\mathrm{out}}(x,y,T)=\sup\{\delta: (x, y+\delta \un_n)\in T)\}
=\sup\{\delta: (\mathbf e^{\alpha x}, e^{ \alpha \delta }\mathbf e^{ \alpha y})\in \mathbf e^{\alpha T)}\} \notag 
\end{equation}
Setting $ \theta =e^{ \alpha \delta}$ we deduce that $\delta = \frac{1}{\alpha}\ln (\theta)$. Therefore:
$$
\mathbb D_{\mathrm{out}}(x,y,T)=\sup \big \{ \frac{1}{\alpha}\ln (\theta): (  \mathbf e^{\alpha x}, \theta \mathbf e^{ \alpha y})\in \mathbf e^{\alpha T}\big\}  =  \frac{1}{\alpha}\ln\Big( \sup \big \{\theta: (  \mathbf e^{\alpha x}, \theta \mathbf e^{ \alpha y})\in \mathbf e^{\alpha T}\big\}\Big),
$$
which completes the proof. $\Box$\\

\subsection{Quantized Non-Parametric Production Models}

In the following we say that a subset $C$ of $\mathbb M_\alpha^n$ is $\mathbb Q_\alpha$-convex if $\sigma_\alpha(C)=\phi_\alpha \circ \mathbf{exp}(C)$ is convex. Along this line we consider a specific class of non-parametric production models. We term it \textit{quantized technology}. The analogue case of constant returns to scale (CRS) in DEA models is expressed as:
 \begin{equation}\notag
 T_{\mathbb Q_\alpha, C}(A)=\Big\{(x,y)\in \Real_+^{m+n}: x\geq \stackrel{\alpha}{\bigoplus_{k\in [\ell]}}t_kx_k,y\leq \stackrel{\alpha}{\bigoplus_{k\in [\ell]}}t_ky_k ,t\in \mathbb M_\alpha^\ell\Big \} 
 \end{equation}
where $A=\{(x_1,y_1),\ldots,(x_{\ell}, y_{\ell})\}\subset \Real_+^{m+n}$ is a set of the observed production vectors. The analogue case of variable returns to scale (VRS) is defined as:
  \begin{equation}\notag
 T_{\mathbb Q_\alpha, V}(A)=\Big\{(x,y)\in \Real_+^{m+n}: x\geq \stackrel{\alpha}{\bigoplus_{k\in [\ell]}}t_kx_k,y\leq \stackrel{\alpha}{\bigoplus_{k\in [\ell]}}t_ky_k, {\stackrel{\alpha}{\bigoplus_{k\in [\ell]}}t_k=0}, \ t\in \mathbb M_\alpha^\ell\Big\}
 \end{equation}
Defining the $\mathbb Q_\alpha$-canonical hull of a finite subset  $A=\{z_1,\ldots,z_\ell\} \subset \Real^d$ as, 
$$
Cc^{\sigma_{\alpha}}(A)=\Big\{ \stackrel{\alpha}{\bigoplus _{k\in[\ell]}}t_{k}\otimes z_{k} :  t\in \mathbb M_\alpha^\ell \Big\},
$$
the following relation arises:
\begin{equation}\notag
T_{\mathbb Q_\alpha, C}(A)=\Big(Cc^{\sigma_\alpha}(A)+K\Big) \cap \Real_{+}^{m+n}
\end{equation}
It is further shown that it satisfies an assumption of Graph-translation homotheticity (T4). The production set $T_{\mathbb Q_\alpha,V}$, being an analogue of VRS, is constructed from the generalized convex hull of the set $A=\{(x_1,y_1),\ldots,(x_{\ell}, y_{\ell})\}\subset \Real_+^{m+n}$:
\begin{equation}\notag
T_{\mathbb Q_\alpha, V}(A)=\Big(Co^{\sigma_\alpha}(A)+K\Big)\cap \Real_{+}^{m+n}
\end{equation}
\bigskip

\begin{center}{\scriptsize 
\unitlength 0.4mm 
\linethickness{0.4pt}
\ifx\plotpoint\undefined\newsavebox{\plotpoint}\fi 
\begin{picture}(275.25,120.25)(0,0)
\put(4,11.5){\vector(0,1){100.75}}
\put(140,11.5){\vector(0,1){100.75}}
\put(4,11.25){\vector(1,0){102.5}}
\put(140,11.25){\vector(1,0){102.5}}
\put(28.5,42.5){\circle*{1.803}}
\put(164.5,42.5){\circle*{1.803}}
\put(66.75,86.5){\circle*{1.803}}
\put(57.75,54.5){\circle*{1.803}}
\put(194.25,54){\circle*{1.803}}
\put(202.75,86.5){\circle*{1.803}}
\put(38,71){\circle*{1.803}}
\put(174,71){\circle*{1.803}}
\put(20.25,26.75){\circle*{1.581}}
\put(156.25,26.75){\circle*{1.581}}
\put(113.25,11){\makebox(0,0)[cc]{$x$}}
\put(249.25,11){\makebox(0,0)[cc]{$x$}}
\put(21.5,44.75){\makebox(0,0)[cc]{$z_1$}}
\put(158,44.25){\makebox(0,0)[cc]{$z_1$}}
\put(66,91.75){\makebox(0,0)[cc]{$z_2$}}
\put(202.5,91.25){\makebox(0,0)[cc]{$z_2$}}
\put(57.75,48){\makebox(0,0)[cc]{$z_2$}}
\put(194.25,47.5){\makebox(0,0)[cc]{$z_2$}}
\put(25.5,23.25){\makebox(0,0)[cc]{$z_5$}}
\put(162,22.75){\makebox(0,0)[cc]{$z_5$}}
\put(32.75,77){\makebox(0,0)[cc]{$z_4$}}
\put(169.25,76.5){\makebox(0,0)[cc]{$z_4$}}
\put(3.75,120.25){\makebox(0,0)[cc]{$y$}}
\put(139,119.75){\makebox(0,0)[cc]{$y$}}
\put(0,7.75){\makebox(0,0)[cc]{$0$}}
\put(136,7.75){\makebox(0,0)[cc]{$0$}}
\put(69.5,0){\makebox(0,0)[cc]
{{{\bf Figure \ref{Quantized}.1:} Quantized Technology  $\alpha>0$.\ \ \ \ \ }}}
\put(205.5,0){\makebox(0,0)[cc]
{{{\bf Figure \ref{Quantized}.2:} Quantized Technology  $\alpha<0$.}}}
\put(66.75,86.75){\line(1,0){33.25}}
\put(202.75,86.75){\line(1,0){33.25}}
\put(20.25,27){\line(0,-1){15.5}}
\put(156.25,27){\line(0,-1){15.5}}
\put(82.25,61.25){\makebox(0,0)[cc]{$T_{\mathbb Q_\alpha,V}$}}
\put(220.5,60.5){\makebox(0,0)[cc]{$T_{\mathbb Q_{\alpha},V}$}}
\qbezier(37.5,71.25)(28.625,64.125)(28.25,43.5)
\qbezier(164.25,43)(173.125,50.125)(173.5,70.75)
\qbezier(28.75,41.75)(20.5,36)(20.25,27.25)
\qbezier(155.75,27.75)(164,33.5)(164.25,42.25)
\qbezier(37.75,71.75)(49.25,71.625)(66.75,87)
\qbezier(202,86.5)(190.5,86.625)(173,71.25)
\end{picture}

}
\end{center}
\bigskip

The quantized CRS and VRS technologies satisfy the following properties.

\begin{prop}Suppose that $A=\{(x_1,y_1),\ldots,(x_{\ell}, y_{\ell})\}\subset \Real_{+ }^{m+n}$. We have the following properties:

 $(a)$ $T_{\mathbb Q_\alpha, C}(A)$ and $T_{\mathbb Q_\alpha, V}(A) $ are closed \emph{(T1)}. 

 $(b)$ $T_{\mathbb Q_\alpha, C}(A)$ and $T_{\mathbb Q_\alpha, V}(A) $ are $\mathbb Q_\alpha$-convex.  
 $(c)$ $T_{\mathbb Q_\alpha, C}(A)$ and $T_{\mathbb Q_\alpha, V}(A) $satisfy a free disposal assumption \emph{(T3)}. 

 $(d)$ $T_{\mathbb Q_\alpha, C}(A)$ is Graph-translation homothetic \emph{(T4)}. 
\end{prop}
\noindent{\bf Proof:} $(a)$ Since $A\subset \Real_{++}^d$ and since the map $\mathbf{ln}: \Real_{++}^d\mapsto \Real^d$ and  $\mathbf{e }: \Real_{}^d\mapsto \Real_{++}^d$ are continuous, it follows that $Co^{\sigma_\alpha}(A)$ and $Cc^{\sigma_\alpha}(A)$ are closed. Since by definition
$T_{\mathbb Q_\alpha, C}(A)=\big(Cc^{\sigma_\alpha}(A)+K\big)\cap \Real_{+}^{m+n}. 
$ and $T_{\mathbb Q_\alpha, V}(A)=\big(Co^{\sigma_\alpha}(A)+K\big)\cap \Real_{+}^{m+n} $ , it follows from \cite{bl11} that $T_{\mathbb Q_\alpha, C}$ and $T_{\mathbb Q_\alpha, V}$ are closed.  $(b)$ The convex cone $K$ and $\Real_+^d$ are $\mathbb Q_\alpha$-convex   and from the expressions of $T_{\mathbb Q_\alpha, C}$ and $T_{\mathbb Q_\alpha, V}$, it follows that $T_{\mathbb Q_\alpha, C}(A)$ and $T_{\mathbb Q_\alpha, V}(A)$ are $\mathbb Q_\alpha$-convex. $(c)$ It is immediate from the expressions of $T_{\mathbb Q_\alpha, C}$ and $T_{\mathbb Q_\alpha, V}$. $(d)$ Suppose that $(x,y)\in T_{\mathbb Q_\alpha, C}$, then by definition,
$$
x\geq \stackrel{\alpha}{\bigoplus_{k\in [\ell]}}t_kx_k \text{ \ and \ } y\leq \stackrel{\alpha}{\bigoplus_{k\in [\ell]}}t_ky_k  \text{ \ for some $ t\in \mathbb M_\alpha^\ell$.}
$$
Equivalently:
$$
x\geq \frac{1}{\alpha}\mathbf{ln}\Big({\sum_{k\in [\ell]}}e^{\alpha t_k}\mathbf e^{\alpha x_k}\Big)\quad \text{ and }\quad y\leq \frac{1}{\alpha}\mathbf{ln}\Big({\sum_{k\in [\ell]}}e^{\alpha t_k}\mathbf e^{\alpha y_k}\Big) 
$$
It follows that for all real numbers $\delta$,
$$
x+\delta \un_m\geq  \delta \un_m +\frac{1}{\alpha}\mathbf{ln}\Big({\sum_{k\in [\ell]}}e^{\alpha t_k}\mathbf e^{\alpha x_k}\Big)\quad \text{ and }\quad y+\delta \un_n \leq \delta \un_n+ \frac{1}{\alpha}\mathbf{ln}\Big({\sum_{k\in [\ell]}}e^{\alpha t_k}\mathbf e^{\alpha y_k}\Big) 
$$
Therefore,
$$
x+\delta \un_m\geq \frac{1}{\alpha } \mathbf{ln}\big( \mathbf e^{\alpha \delta \un_m}\big)+\frac{1}{\alpha}\mathbf{ln}\Big({\sum_{k\in [\ell]}}e^{\alpha t_k}\mathbf e^{\alpha x_k}\Big)\quad \text{ and }\quad y +\delta \un_n \leq \frac{1}{\alpha } \mathbf{ln}\big( \mathbf e^{\alpha \delta  \un_n}\big)
 + \frac{1}{\alpha}\mathbf{ln}\Big({\sum_{k\in [\ell]}}e^{\alpha t_k}\mathbf e^{\alpha y_k}\Big)
 $$
Finally,
$$
x+\delta \un_m\geq \frac{1}{\alpha}\mathbf{ln}\Big({\sum_{k\in [\ell]}}e^{\alpha (t_k+\delta)}\mathbf e^{\alpha x_k}\Big)\quad \text{ and }\quad y +\delta \un_n \leq   \frac{1}{\alpha}\mathbf{ln}\Big({\sum_{k\in [\ell]}}e^{\alpha (t_k+\delta)}\mathbf e^{\alpha y_k}\Big)
$$
This implies that $(x+\delta \un_m, y +\delta \un_n )\in T_{\mathbb Q_\alpha, C}$. Consequently, $T_{\mathbb Q_\alpha, C}$ satisfies Graph-translation homothetic (T4). $\Box$\\
 
\begin{cor}Suppose that $A\subset \Real^d$, then the following properties arise:

$(a)$ $\mathbb D_{\mathrm{in}}\big((x,y,T_{\mathbb Q_\alpha, V}(A)\big)=-\frac{1}{\alpha}\ln\big(E_{\mathrm{in}}(\mathbf e^{\alpha x}, \mathbf e^{\alpha y }, T_V(\mathbf{e}^{\alpha A}))\big)$
 
 $(b)$ $ \mathbb D_{\mathrm{out}}\big((x,y,T_{\mathbb Q_\alpha, V}(A)\big)=-\frac{1}{\alpha}\ln\big(E_{\mathrm{out}}(\mathbf e^{\alpha x}, \mathbf e^{\alpha y }, T_V(\mathbf{e}^{\alpha A}))\big)$
 
  $(c)$ $ \mathbb D_{\mathrm{in}}\big((x,y,T_{\mathbb Q_\alpha, C}(A)\big)=-\frac{1}{\alpha}\ln\big(E_{\mathrm{in}}(\mathbf e^{\alpha x}, \mathbf e^{\alpha y }, T_C(\mathbf{e}^{\alpha A}))\big)$
  
 $(d)$ $\mathbb D_{\mathrm{out}}\big((x,y,T_{\mathbb Q_\alpha,C }(A)\big)=-\frac{1}{\alpha}\ln\big(E_{\mathrm{out}}(\mathbf e^{\alpha x}, \mathbf e^{\alpha y }, T_C(\mathbf{e}^{\alpha A}))\big)$
 \end{cor}
{\bf Proof:} We have the relation $T_V(\mathbf{e}^{\alpha A}) =\mathbf e^{T_{\mathbb Q_\alpha, V}(A)}$ and 
$T_C(\mathbf{e}^{\alpha A}) =\mathbf e^{T_{\mathbb Q_\alpha, C}(A)}$. The result immediately follows from Proposition \ref{Distance}. $\Box$\\

In the following we show how to convert the problem of measuring efficiency into a linear programming problem: 
\begin{align}\notag
\begin{array}{llllll}
\mathbb D_{\mathrm{in}}( x,y, T_{\mathbb Q_\alpha, V})  = & \max \delta \\
\text{subject to} \qquad & x-\delta \un_m  \geq \frac{1}{\alpha }\mathbf{ln}\Big(\sum\limits_{k \in [\ell]} e^{\alpha t_k}  {\mathbf e^{\alpha x_k} }\Big)  \\
& y \leq \frac{1}{\alpha }\mathbf{ln}\Big(\sum\limits_{k \in [\ell]}  e^{\alpha t_k}  {\mathbf e^{\alpha y_k} }\Big)  \\
& \frac{1}{\alpha }\ln\Big(\sum\limits_{k \in [\ell]}  e^{\alpha t_k}  \Big)=0,\  {t\in \mathbb M_{\alpha}^\ell}
\end{array}
\end{align}
Applying the transformation $z\mapsto \mathbf{e}^{\alpha z} $ to both sides of each each equation yields:
\begin{align}\notag
\begin{array}{llllll}
\mathbb D_{\mathrm{in}}( x,y, T_{\mathbb Q_\alpha, V})  = & \max \delta \\
\text{subject to:} \qquad & \mathbf e^{-\alpha \delta} \mathbf{e}^{\alpha x}   \geq  \sum\limits_{k \in [\ell]} e^{\alpha t_k} {\mathbf e^{\alpha x_k} }   \\
&  \qquad \mathbf e^{\alpha y} \leq  \sum\limits_{k \in [\ell]}  e^{\alpha t_k}  {\mathbf e^{\alpha y_k} }   \\
&  \sum\limits_{k \in [\ell]}  e^{\alpha t_k}   =1,\  {t\in \mathbb M_{\alpha}^\ell}
\end{array}
\end{align}
Setting $s _k= e^{\alpha t_k}$ for all $k\in [\ell]$ and $\lambda =  \mathbf{e}^{-\alpha \delta}$, the following linear program can be deduced from the previous one:
\begin{align}\notag
\begin{array}{llllll}
\mathbb D_{\mathrm{in}}( x,y, T_{\mathbb Q_\alpha, V})  = & \min \lambda \\
\text{subject to:} \qquad & \lambda \mathbf{e}^{\alpha x}   \geq  \sum\limits_{k \in [\ell]} s_k {\mathbf e^{\alpha x_k} }  \\
& \;\mathbf e^{\alpha y} \leq  \sum\limits_{k \in [\ell]}  s_k {\mathbf e^{\alpha y_k} }  \\
&  \sum\limits_{k \in [\ell]} s_k=1,\   {s \geq 0}
\end{array}
\end{align}

\subsection{Duality}

This section is devoted to the dual properties of distance functions, \textit{i.e.} the link between distance functions with revenue and cost functions. 

\noindent The \textit{quantized inner product} $\langle \cdot, \cdot \rangle _{\mathbb Q_\alpha}: \mathbb M_{\alpha}^d\times \mathbb M_{\alpha}^d \mapsto \mathbb M_\alpha$ is defined as:
\begin{equation}\notag
\langle v, z \rangle _{\mathbb Q_\alpha}=\frac{1}{\alpha}\ln \big(\sum_{k\in [d]}e^{\alpha v_k} {e^{\alpha z_k}}\big)
\end{equation}

\noindent The \textit{quantized cost function}  {$C_{\mathbb Q_\alpha}: \mathbb M_{\alpha}^m\times \mathbb M_{\alpha}^m \times \mathcal T$} is: 
\begin{equation}\notag
C_{\mathbb Q_\alpha}(w, y, T)=\inf \big\{\langle w, x \rangle _{\mathbb Q_\alpha}: x\in L(y)\big\}
\end{equation}
with $C_{\mathbb Q_\alpha}(w, y, T)=+\infty$ if $L(y)=\emptyset$, {and with $C(w, y, T)$ the cost function based on the usual inner product.} 
 
\noindent The \textit{quantized revenue function} {$R_{\mathbb Q_\alpha}: \mathbb M_{\alpha}^n\times \mathbb M_{\alpha}^n \times \mathcal T$} is: 
 \begin{equation}\notag
R_{\mathbb Q_\alpha}(p, x, T)=\sup \big\{\langle p, y \rangle _{\mathbb Q_\alpha}: y\in P(x)\big\}
\end{equation}
with $R_{\mathbb Q_\alpha}(p,x,T)=-\infty$ if $P(x)=\emptyset$,  {and with $R(p, x, T)$ the revenue function based on the usual inner product.} 

\begin{prop} Let $T\in \mathcal T$ and suppose that $(x,y)\in T$. 
\newline $(i)$ If $L(y)$ is $\mathbb Q_\alpha$-convex, then 
\begin{equation*}
\mathbb D_{\mathrm{in}}(x,y, T)=\inf_{w\in \mathbb M_\alpha^m}  \big\{\langle w, x \rangle _{\mathbb Q_\alpha}-C_{\mathbb Q_\alpha}(w, y, T):{\langle w, 0 \rangle _{\mathbb Q_\alpha}=0} \} 
\end{equation*}
\newline $(ii)$ If $P(x)$ is $\mathbb Q_\alpha$-convex, then 
\begin{equation*}\mathbb D_{\mathrm{out}}(x,y, T)=\inf_{p\in \mathbb M_\alpha^n}  \big\{R_{\mathbb Q_\alpha}{(p, x, T)} - \langle p, y \rangle_{\mathbb Q_\alpha}: {\langle p, 0 \rangle _{\mathbb Q_\alpha}=0}\} 
\end{equation*}
\end{prop}
\noindent{\bf Proof:} $(i)$ We have shown that,
 $$
 \mathbb D_{\mathrm{in}}(x,y,T)=-\frac{1}{\alpha}\ln\Big(E_{\mathrm{in}}(\mathbf e^{\alpha x}, \mathbf e^{\alpha y }, \mathbf{e}^{\alpha T})\Big)
 $$
 However, we have:
 $$
 E_{\mathrm{in}}(\mathbf e^{\alpha x}, \mathbf e^{\alpha y }, \mathbf{e}^{\alpha T})=\sup_{{v\in \mathbb M_\alpha^m}}\Big\{\frac{C(v, \mathbf e^{\alpha y}, \mathbf e^{\alpha T})}{\langle v, \mathbf e^{\alpha x}\rangle }:\langle v ,\un_m\rangle =1\Big\}
 $$
 Setting $v=\mathbf e^{\alpha w}$, we obtain 
  $$
  E_{\mathrm{in}}(\mathbf e^{\alpha x}, \mathbf e^{\alpha y }, \mathbf{e}^{\alpha T})=\sup_{w\in \mathbb M_{\alpha}^m}\Big\{\frac{C(\mathbf e^{\alpha w}, \mathbf e^{\alpha y}, \mathbf e^{\alpha T})}{\langle \mathbf e^{\alpha w}, \mathbf e^{\alpha x}\rangle }:\langle {\mathbf e^{\alpha w} },\un_m\rangle =1\Big\}
$$
Therefore,
 $$
  \mathbb D_{\mathrm{in}}(x,y,T)=-\frac{1}{\alpha}\ln\Big(\sup_{w\in \mathbb M_{\alpha}^m}\Big\{\frac{C(\mathbf e^{\alpha w}, \mathbf e^{\alpha y}, \mathbf e^{\alpha T})}{\langle \mathbf e^{\alpha w}, \mathbf e^{\alpha x}\rangle }:\langle {\mathbf e^{\alpha w} } ,\un_m\rangle =1\Big\}\Big)
$$
It can be deduced that,
$$
\mathbb D_{\mathrm{in}}(x,y,T)= \inf_{{w\in \mathbb M_{\alpha}^m}}\Big\{\frac{1}{\alpha}\ln {\langle \mathbf e^{\alpha w}, \mathbf e^{\alpha x}\rangle }-\frac{1}{\alpha}\ln \Big( C(\mathbf e^{\alpha w}, \mathbf e^{\alpha y}, \mathbf e^{\alpha T})\Big):  \frac{1}{\alpha}\ln  \langle \mathbf e^{\alpha w} ,\mathbf e^0\rangle  =0\Big\} 
$$
However:
  $$\frac{1}{\alpha}\ln \Big( C(\mathbf e^{\alpha w}, \mathbf e^{\alpha y}, \mathbf e^{\alpha T})\Big)=\frac{1}{\alpha}\ln \Big( \inf \big \{\langle \mathbf e^{\alpha w}, u \rangle : (u, \mathbf e^{\alpha y})\in \mathbf e^{\alpha T}\big \}\Big)$$
Setting $u=\mathbf e^{\alpha x}$, yields 
   $$\frac{1}{\alpha}\ln \Big( C(\mathbf e^{\alpha w}, \mathbf e^{\alpha y}, \mathbf e^{\alpha T})\Big)=\inf \Big( \frac{1}{\alpha}\ln \big \{\langle \mathbf e^{\alpha w}, \mathbf e^{\alpha x}\rangle : (x, y)\in T \big \}\Big)=C_{\mathbb Q_\alpha}(w, y, T)
   $$
   Consequently,
   \begin{equation*}\mathbb D_{\mathrm{in}}(x,y, T)=\inf_{w\in \mathbb M_{\alpha}^m}  \big\{\langle w, x \rangle _{\mathbb Q_\alpha}-C_{\mathbb Q_\alpha}(w, y, T):\langle w, 0\rangle _{\mathbb Q_\alpha}=0\} \end{equation*}
$(ii)$ We have established that:
$$
\mathbb D_{\mathrm{out}}(x,y,T)=\frac{1}{\alpha}\ln\Big(E_{\mathrm{out}}(\mathbf e^{\alpha x}, \mathbf e^{\alpha y}, \mathbf{e}^{\alpha T})\Big)
$$
However, we have:
$$
E_{\mathrm{out}}(\mathbf e^{\alpha x}, \mathbf e^{\alpha y }, \mathbf{e}^{\alpha T})=\inf_{{r\in \mathbb M_\alpha^n}}\Big\{\frac{R(r,{\mathbf e^{\alpha x}}, \mathbf e^{\alpha T})}{\langle r, \mathbf e^{\alpha y}\rangle }:\langle r ,\un_n\rangle =1\Big\}
$$
Setting $r=\mathbf e^{\alpha p}$, we obtain:
$$
E_{\mathrm{out}}(\mathbf e^{\alpha x}, \mathbf e^{\alpha y }, \mathbf{e}^{\alpha T})=\inf_{p\in \mathbb M_\alpha^n}\Big\{\frac{R(\mathbf e^{\alpha p}, \mathbf e^{\alpha x}, \mathbf e^{\alpha T})}{\langle \mathbf e^{\alpha p}, \mathbf e^{\alpha y}\rangle }:\langle {\mathbf e^{\alpha p}},\un_n\rangle = 1\Big\}
$$
Therefore,
$$
  \mathbb D_{\mathrm{out}}(x,y,T)=-\frac{1}{\alpha}\ln\Big(\sup_{p\in \mathbb M_\alpha^m}\Big\{\frac{R(\mathbf e^{\alpha p}, \mathbf e^{\alpha x}, \mathbf e^{\alpha T})}{\langle \mathbf e^{\alpha p}, \mathbf e^{\alpha y}\rangle }:\langle {\mathbf e^{\alpha p}} ,\un_n\rangle =1\Big\}\Big)
$$
We deduce that,
$$
\mathbb D_{\mathrm{out}}(x,y,T)= \inf_{p\in \mathbb M_\alpha^n}\Big\{ { \frac{1}{\alpha}\ln \Big( R(\mathbf e^{\alpha p}}, \mathbf e^{\alpha x}, \mathbf e^{\alpha T})\Big){-\frac{1}{\alpha}\ln} \langle \mathbf e^{\alpha p}, \mathbf e^{\alpha y}\rangle :  \frac{1}{\alpha}\ln  \langle {\mathbf e^{\alpha p}} ,\mathbf e^0\rangle  =0\Big\}
$$
However:
$$
\frac{1}{\alpha}\ln \Big( R(\mathbf e^{\alpha p}, \mathbf e^{\alpha x}, \mathbf e^{\alpha T})\Big)=\frac{1}{\alpha}\ln \Big( \inf \big \{\langle {\mathbf e^{\alpha p}}, h \rangle : (\mathbf e^{\alpha x},h)\in \mathbf e^{\alpha T}\big \}\Big)
$$
Setting $h=\mathbf e^{\alpha y}$, yields 
$$
\frac{1}{\alpha}\ln \Big( R(\mathbf e^{\alpha p}, \mathbf e^{\alpha x}, \mathbf e^{\alpha T})\Big)=\inf \Big \{  {\frac{1}{\alpha}\ln } \langle \mathbf e^{\alpha p}, \mathbf e^{\alpha y}\rangle : (x, y)\in T \Big \}=R_{\mathbb Q_\alpha}(p, x, T)
$$
Finally:
\begin{equation*}
\mathbb D_{\mathrm{out}}(x,y, T)=\inf_{p\in \mathbb M_\alpha^n}  \big\{R_{\mathbb Q_\alpha}(p, x, T)-\langle p, y \rangle _{\mathbb Q_\alpha} :\langle {p}, 0\rangle _{\mathbb Q_\alpha}=0\} \quad \Box
\end{equation*}

\section{Dequantization and Tropical Limits: A Discrete Model}\label{Dequantized}
   
This section proposes computational formula for measuring efficiency with tropical limits of production technologies (Subsection 1). Closed formula for distance functions are proposed in the case of non-parametric tropical technologies (Subsection 2). Then, it is shown that these formula apply for a large class of discrete production models (Subsection 3) such as FDH ones (Subsection 4).  

\subsection{Tropical Limit of Production Technologies}

Given a subset $C$ of $\Real^d$, $ Hh(C)$ denotes the set homogeneously spanned from $C$, equivalently $ Hh(C)=\{\lambda
 {v}:  {v}\in C,\lambda \in \Real\}$. The translation homothetic spanned set is defined a $Th(C)=\{v+ \delta \un_d: v\in C\}=\{ \delta \un_d \otimes v: v\in C\}$.
Along this line, let us define,
\begin{equation}\label{TarasHull}
Cc^{\sigma_\infty}(A)=Co^{\sigma_\infty}(Th(A))\quad \text{ and }\quad Cc^{\sigma_{-\infty}}(A)=Co^{\sigma_{-\infty}}(Th(A))
\end{equation}
Notice that we have the relations, for all $\alpha\in \Real$, 
\begin{equation}\label{DefTrans}
 {Cc^{\sigma_\alpha}(A)}=Co^{\sigma_\alpha}\big (Hh(A)\big)
\end{equation}
and
\begin{equation}\label{LonTrans}
Co^{\sigma_\alpha}(A)=\mathbf{ln}\Big(Co^{\alpha}(\mathbf{exp}(A))\Big)
\end{equation}
It follows that if $A\in\mathbb M_\alpha^{d}$, then:
 
\begin{equation}\label{ConTrans}
Cc^{\sigma_\alpha}(A)=\mathbf{ln}\Big( {Cc^{ \alpha}}(\mathbf{exp}(A))\Big)
\end{equation}
Moreover, we have the relation,
\begin{equation}\label{TransHom}
 \mathbf{exp}\big(Th(A)\big)=Hh\big (\mathbf{ exp}(A)\big)
\end{equation}
Paralleling our earlier definitions let us define the following tropical technologies, see \cite{abs17}:
\begin{equation}\notag
   T_{\mathbb Q_\infty, V}(A)=\Big(Co^{\sigma_\infty}(A)+K\Big)\cap \Real_+^d
\end{equation}
This model corresponds to the upper-dequantization of the quantized  {VRS production model} taking the limit when $\alpha\longrightarrow \infty$. 
 This set can be equivalently defined as:
\begin{equation}\notag
T_{\mathbb Q_{\infty},V} (A)  =\Big\{(x,y)\in \mathbb{R}_{+}^{d}:
 x\geq \stackrel{+\infty}{\bigoplus_{k\in [\ell]}}(t_{k}\otimes x_{k}),
 \,y\leq \stackrel{+\infty}{\bigoplus_{k\in [\ell]}}(t_{k}\otimes y_{k}), \max_{k\in [\ell]}t_{k}=0,
  t \in \mathbb M_\infty^\ell \Big\}
 \end{equation}
 {On the other hand, the tropical graph translation homothetic model corresponds to the upper-dequantization of the quantized CRS model taking the limit when $\alpha\longrightarrow \infty$:}

\begin{equation}\notag
    {T_{\mathbb Q_{\infty},C}(A)}=\Big(Cc^{\sigma_\infty}(A)+K\Big)\cap \Real_+^d
\end{equation}
Equivalently,
\begin{equation}\notag
T_{\mathbb Q_\infty, C} (A)  =\Big\{(x,y)\in \mathbb{R}_{+}^{d}:
 x\geq \stackrel{+\infty}{\bigoplus_{k\in [\ell]}}(t_{k}\otimes x_{k}),
 \,y\leq \stackrel{+\infty}{\bigoplus_{k\in [\ell]}}(t_{k}\otimes y_{k}),  
  t \in \mathbb M_\infty^\ell \Big\}
 \end{equation}
It is therefore immediate to consider the lower-dequantization  {of the quantized VRS model} that is obtained taking the limit when $\alpha\longrightarrow -\infty$. It is defined as:
  \begin{equation}\notag
   T_{\mathbb Q_{-\infty},V}(A)=\Big(Co^{\sigma_{-\infty}}(A)+K\Big)\cap \Real_+^d
   \end{equation}
Equivalently, this model can be formulated as follows:
\begin{equation}\notag
T_{\mathbb Q_{-\infty}, V}(A)   =\Big\{(x,y)\in \mathbb{R}_{+}^{d}:
 x\geq \stackrel{-\infty}{\bigoplus_{k\in [\ell]}}(t_{k}\otimes x_{k}),
 \,y\leq \stackrel{-\infty}{\bigoplus_{k\in [\ell]}}(t_{k}\otimes y_{k}), \min_{k\in [\ell]}t_{k}=0,
  t \in \mathbb M_{-\infty}^\ell \Big\}
 \end{equation}
  {The tropical graph translation homothetic model, corresponding to the lower-dequantization of the quantized CRS model, is defined as:}
 
 \begin{equation}\notag
   T_{\mathbb Q_{-\infty},C}(A)=\Big(Cc^{\sigma_{-\infty}}(A)+K\Big)\cap \Real_+^d
   \end{equation}
Equivalently, 
\begin{equation}\notag
T_{\mathbb Q_{-\infty}, C}(A ) =\Big\{(x,y)\in \mathbb{R}_{+}^{d}:
 x\geq \stackrel{-\infty}{\bigoplus_{k\in [\ell]}}(t_{k}\otimes x_{k}),
 \,y\leq \stackrel{-\infty}{\bigoplus_{k\in [\ell]}}(t_{k}\otimes y_{k}), 
  t \in \mathbb M_{-\infty}^\ell \Big\}
\end{equation}
\medskip
 
\begin{center}
{\scriptsize \unitlength 0.4mm 
\linethickness{0.4pt}
\ifx\plotpoint\undefined\newsavebox{\plotpoint}\fi 
\begin{picture}(275.25,120.25)(0,0)
\put(4,11.5){\vector(0,1){100.75}}
\put(140,11.5){\vector(0,1){100.75}}
\put(4,11.25){\vector(1,0){102.5}}
\put(140,11.25){\vector(1,0){102.5}}
\put(28.5,42.5){\circle*{1.803}}
\put(164.5,42.5){\circle*{1.803}}
\put(66.75,86.5){\circle*{1.803}}
\put(57.75,54.5){\circle*{1.803}}
\put(194.25,54){\circle*{1.803}}
\put(202.75,86.5){\circle*{1.803}}
\put(38,71){\circle*{1.803}}
\put(174,71){\circle*{1.803}}
\put(20.25,26.75){\circle*{1.581}}
\put(156.25,26.75){\circle*{1.581}}
\put(113.25,11){\makebox(0,0)[cc]{$x$}}
\put(249.25,11){\makebox(0,0)[cc]{$x$}}
\put(21.5,44.75){\makebox(0,0)[cc]{$z_1$}}
\put(158,44.25){\makebox(0,0)[cc]{$z_1$}}
\put(66,91.75){\makebox(0,0)[cc]{$z_2$}}
\put(202.5,91.25){\makebox(0,0)[cc]{$z_2$}}
\put(57.75,48){\makebox(0,0)[cc]{$z_2$}}
\put(194.25,47.5){\makebox(0,0)[cc]{$z_2$}}
\put(25.5,23.25){\makebox(0,0)[cc]{$z_5$}}
\put(162,22.75){\makebox(0,0)[cc]{$z_5$}}
\put(32.75,77){\makebox(0,0)[cc]{$z_4$}}
\put(169.25,76.5){\makebox(0,0)[cc]{$z_4$}}
\put(3.75,120.25){\makebox(0,0)[cc]{$y$}}
\put(139,119.75){\makebox(0,0)[cc]{$y$}}
\put(0,7.75){\makebox(0,0)[cc]{$0$}}
\put(136,7.75){\makebox(0,0)[cc]{$0$}}
\put(69.5,0){\makebox(0,0)[cc]
{{{\bf Figure \ref{Dequantized}.1:} Max-Plus  Production Set.}}}
\put(205.5,0){\makebox(0,0)[cc]{{{\bf Figure \ref{Dequantized}.2:} Min-Plus  Production Set.}}}
\put(28.75,63.75){\line(0,-1){21.5}}
\put(173.75,51){\line(0,1){21.5}}
\put(38,71.25){\line(1,0){11.5}}
\put(202,86.5){\line(-1,0){11.5}}
\multiput(37.25,71.25)(-.095108696,-.08423913){92}{\line(-1,0){.095108696}}
\multiput(165.25,43.5)(.095108696,.08423913){92}{\line(1,0){.095108696}}
\multiput(28.5,42.25)(-.09638554,-.08433735){83}{\line(-1,0){.09638554}}
\multiput(156.25,27.5)(.09638554,.08433735){83}{\line(1,0){.09638554}}
\put(20.5,35.25){\line(0,-1){8.75}}
\put(164.25,34.5){\line(0,1){8.75}}
\multiput(49.5,71.5)(.095303867,.084254144){181}{\line(1,0){.095303867}}
\multiput(190.5,86.25)(-.095303867,-.084254144){181}{\line(-1,0){.095303867}}
\put(66.75,86.75){\line(1,0){33.25}}
\put(202.75,86.75){\line(1,0){33.25}}
\put(20.25,27){\line(0,-1){15.5}}
\put(156.25,27){\line(0,-1){15.5}}
\put(3.75,11.75){\vector(1,1){9.25}}
\put(140,11){\vector(1,1){9.25}}
\put(82.25,61.25){\makebox(0,0)[cc]{$T_{\mathbb Q_\infty,V} $}}
\put(220.5,60.5){\makebox(0,0)[cc]{$T_{\mathbb Q_{-\infty},V}$}}
\end{picture}}
\end{center}
\bigskip

In order to establish the next results, we need the two following intermediate results established in \cite{abm17}. 
 
\begin{prop}\label{conlim} Let
$\{C_n\}_{n\in\mathbb{N}}$, be a sequence of compact sets of points of $\Real^d$ which converges, in the
Painlev\'e-Kuratowski sense, to a set $C$ of $\Real^d$, that is $\Lim_{n\longrightarrow \infty}C_n=C$. Then, we have:
$$
\Lim_{n\longrightarrow \infty}Hh(C_n)=Hh(C)
$$
\end{prop}

\begin{prop} \label{interlim} Let $\{C_n\}_{n\in\mathbb{N}}$ be a sequence of compact sets of points of $\Real_+^d$ which converges, in the
Painlev\'e-Kuratowski sense, to a set $C$ of $\Real_+^d$, that is $\Lim_{n\longrightarrow \infty}C_n=C$. Let $K=\Real_+^{m}\times (-\Real_+^{n})$ with $m+n=d$.  Then, we have:
$$
\Lim_{n\longrightarrow \infty}\big[(C_n+K)\cap \Real_+^d\big]=(C+K)\cap \Real_+^d
$$

\end{prop}

  {The tropical production sets, being either upper-dequantized or lower-dequantized, can be expressed in terms of the limit of the quantized non-parametric production sets (VRS or CRS).}
 
 \begin{prop} Suppose that $A=\{(x_1,y_1),\ldots,(x_\ell,y_\ell)\}$ is a finite subset of $\Real^d$. Suppose that $\{\alpha_k\}_{k\in \mathbb N}$ is a sequence of real numbers such that $\lim_{k\longrightarrow \infty}\alpha_k=+\infty$. Then,
 $$
 \Lim_{k\longrightarrow \infty}T_{\mathbb Q_{\alpha_k}, V}(A)=T_{\mathbb Q_\infty, V}(A)  
 $$
 and
 $$
 \Lim_{k\longrightarrow \infty}T_{\mathbb Q_{\alpha_k}, C}(A)=T_{\mathbb Q_\infty, C}(A)
 $$
Suppose that $\{\alpha_k\}_{k\in \mathbb N}$ is a sequence of real numbers such that $\lim_{k\longrightarrow \infty}\alpha_k=-\infty$. Then,
$$
\Lim_{k\longrightarrow \infty}T_{\mathbb Q_{\alpha_k}, V}(A)=T_{\mathbb Q_{-\infty}, V}(A)  
$$
and
$$
\Lim_{k\longrightarrow \infty}T_{\mathbb Q_{\alpha_k}, C}(A)=T_{\mathbb Q_{-\infty}, C}(A) 
$$
\end{prop}

\noindent{\bf Proof:} To prove the results we use equations \eqref{TarasHull}, \eqref{DefTrans}, \eqref{LonTrans}, \eqref{ConTrans} and \eqref{TransHom}. 
In the VRS case the result is an immediate consequence of Proposition \ref{LimTRopConv} and Proposition \ref{interlim}. Therefore, we have:
 $$
\Lim_{k\longrightarrow \infty}T_{\mathbb Q_{\alpha_k}, V}(A)=T_{\mathbb Q_\infty, V}(A)  
$$
To prove the translation homothetic case we use the fact that for all $k$:
 
\begin{equation*} 
Cc^{\alpha_k}(\mathbf e^A)=Co^{\alpha_k}\big (Hh(\mathbf e^A)\big)
\end{equation*}
Moreover,
 \begin{equation*} 
Cc^{\sigma_{\alpha_k}}(A)=\mathbf{ln}\Big(Cc^{\alpha_k}(\mathbf{exp}(A))\Big)
\end{equation*}
Consequently, since $\mathbf{ln}$ is continuous on $\Real_{++}^d$, using Proposition \ref{conlim}: 
 \begin{align*} 
\Lim_{k\longrightarrow \infty}Cc^{\sigma_{\alpha_k}}(A)&=\Lim_{k\longrightarrow \infty}\mathbf{ln}\Big(Cc^{\alpha_k}(\mathbf e^A )\Big)\\
&=\mathbf{ln}\Lim_{k\longrightarrow \infty}\Big(Cc^{\alpha_k}(\mathbf e^A )\Big)\\&=\mathbf{ln}\Lim_{k\longrightarrow \infty}\Big(Co^{\alpha_k}\big (Hh(\mathbf e^A)\big)\Big)\\&=\mathbf{ln} \Big(Co^{\infty}\big (Hh(\mathbf e^A)\big)\Big)\\&=\mathbf{ln} \Big(Co^{\infty}\big  (\mathbf{exp}(Th(A))\big)\Big)=Cc^{\sigma_\infty}(A)
\end{align*}
Therefore, we deduce that: 
$$
\Lim_{k\longrightarrow \infty}T_{\mathbb Q_{\alpha_k}, C}(A)=T_{\mathbb Q_\infty, C}(A) 
$$
 
The proof of the case where $\lim_{k\longrightarrow \infty}\alpha_k=-\infty$ is obtained using similar arguments. $\Box$. \\

 \begin{center}
 {\scriptsize 
 
\unitlength 0.4mm 
\linethickness{0.4pt}
\ifx\plotpoint\undefined\newsavebox{\plotpoint}\fi 
\begin{picture}(275.25,120.25)(0,0)
\put(4,11.5){\vector(0,1){100.75}}
\put(140,11.5){\vector(0,1){100.75}}
\put(4,11.25){\vector(1,0){102.5}}
\put(140,11.25){\vector(1,0){102.5}}
\put(28.5,42.5){\circle*{1.803}}
\put(164.5,42.5){\circle*{1.803}}
\put(66.75,86.5){\circle*{1.803}}
\put(57.75,54.5){\circle*{1.803}}
\put(194.25,54){\circle*{1.803}}
\put(202.75,86.5){\circle*{1.803}}
\put(38,71){\circle*{1.803}}
\put(174,71){\circle*{1.803}}
\put(20.25,26.75){\circle*{1.581}}
\put(156.25,26.75){\circle*{1.581}}
\put(113.25,11){\makebox(0,0)[cc]{$x$}}
\put(249.25,11){\makebox(0,0)[cc]{$x$}}
\put(21.5,44.75){\makebox(0,0)[cc]{$z_1$}}
\put(158,44.25){\makebox(0,0)[cc]{$z_1$}}
\put(66,91.75){\makebox(0,0)[cc]{$z_2$}}
\put(202.5,91.25){\makebox(0,0)[cc]{$z_2$}}
\put(57.75,48){\makebox(0,0)[cc]{$z_2$}}
\put(194.25,47.5){\makebox(0,0)[cc]{$z_2$}}
\put(25.5,23.25){\makebox(0,0)[cc]{$z_5$}}
\put(162,22.75){\makebox(0,0)[cc]{$z_5$}}
\put(32.75,77){\makebox(0,0)[cc]{$z_4$}}
\put(169.25,76.5){\makebox(0,0)[cc]{$z_4$}}
\put(3.75,120.25){\makebox(0,0)[cc]{$y$}}
\put(139,119.75){\makebox(0,0)[cc]{$y$}}
\put(0,7.75){\makebox(0,0)[cc]{$0$}}
\put(136,7.75){\makebox(0,0)[cc]{$0$}}
\put(69.5,0){\makebox(0,0)[cc]
{{{\bf Figure \ref{Dequantized}.3:} Limit Set $\alpha=+\infty$.}}}
\put(205.5,0){\makebox(0,0)[cc]
{{{\bf Figure \ref{Dequantized}.4:} Limit Set $\alpha=-\infty$.}}}
\put(28.75,63.75){\line(0,-1){21.5}}
\put(173.75,51){\line(0,1){21.5}}
\put(38,71.25){\line(1,0){11.5}}
\put(202,86.5){\line(-1,0){11.5}}
\multiput(37.25,71.25)(-.095108696,-.08423913){92}{\line(-1,0){.095108696}}
\multiput(165.25,43.5)(.095108696,.08423913){92}{\line(1,0){.095108696}}
\multiput(28.5,42.25)(-.09638554,-.08433735){83}{\line(-1,0){.09638554}}
\multiput(156.25,27.5)(.09638554,.08433735){83}{\line(1,0){.09638554}}
\put(20.5,35.25){\line(0,-1){8.75}}
\put(164.25,34.5){\line(0,1){8.75}}
\multiput(49.5,71.5)(.095303867,.084254144){181}{\line(1,0){.095303867}}
\multiput(190.5,86.25)(-.095303867,-.084254144){181}{\line(-1,0){.095303867}}
\put(66.75,86.75){\line(1,0){33.25}}
\put(202.75,86.75){\line(1,0){33.25}}
\put(20.25,27){\line(0,-1){15.5}}
\put(156.25,27){\line(0,-1){15.5}}
\put(3.75,11.75){\vector(1,1){9.25}}
\put(140,11){\vector(1,1){9.25}}
\put(82.25,61.25){\makebox(0,0)[cc]{$T_{\mathbb Q_\infty,V}$}}
\put(220.5,60.5){\makebox(0,0)[cc]{$T_{\mathbb Q_{-\infty},V}$}}
\qbezier(37.5,71.25)(28.625,64.125)(28.25,43.5)
\qbezier(164.25,43)(173.125,50.125)(173.5,70.75)
\qbezier(28.75,41.75)(20.5,36)(20.25,27.25)
\qbezier(155.75,28.25)(164,34)(164.25,42.75)
\qbezier(37.75,71.75)(49.25,71.625)(66.75,87)
\qbezier(202,86.5)(190.5,86.625)(173,71.25)
\end{picture}
}
\end{center}
\medskip

\subsection{Distance Functions and Non-Parametric Models}
 
This section extends some results established in \cite{abs17} to the case of Min-Plus technologies. For this purpose, some dual relationships between Max-Plus and Min-Plus models are proposed. 
 \begin{prop} \label{GaugeMaxPolyVRS}
Let $A=\big\{( x_k,y_k):k\in [\ell]\big\}\subset \mathbb M_{ \infty}^{d}$,  and for all $ \bar k\in [\ell]$ let $\beta_{\bar k,k}=\min\limits_{i=1,\ldots,m}\{x_{\bar k,i}-x_{k,i} \}.$

$(a)$  If
 $ T=T_{\mathbb Q_{ \infty}, V}$ then:
\begin{align}\notag
&\mathbb D_{\mathrm{in}}\big(x_{\bar k},y_{\bar k}, T_{\mathbb Q_{ \infty}, V} (A) \big) =\min\Big\{\min_{j=1,...,n}
\max_{\substack{k\in [\ell] \\ y_{\bar k, j}  \leq y_{k,j} }}\big\{
-y_{\bar k, j}^{}+y_{k,j}^{} + \beta_{\bar k,k}\big\},
\max_{k\in [\ell]}\beta_{\bar k,k}\Big\} \\
&\mathbb D_{\mathrm{out}}\big(x_{\bar k},y_{\bar k}, T_{\mathbb Q_{ \infty},V} (A)\big)=
\min_{j=1,...,n}\max_{\substack{k\in [\ell]}}\left\{y_{k_j} -{y_{\bar
k,j} }+\min \big\{\beta_{\bar k,k}, 0\big\}\right\} \notag
\end{align}

$(b)$ If
 $ T=T_{\mathbb Q_{ \infty}, C}$ then: {\small $$  \mathbb D_{\mathrm{in}}\big (x_{\bar k},y_{\bar
k}, T_{\mathbb Q_{ \infty},C} (A)\big)= \mathbb D_{\mathrm{out}}\big(x_{\bar k},y_{\bar k}, T_{\mathbb Q_{ \infty}, C} (A)\big)
= \min_{j=1,...,n} \max_{k\in [\ell]}\big\{ -y_{\bar k,j}^{}+y_{k,j}^{} + \beta_{\bar
k,k}\big\}$$}
\end{prop}

We now adopt the following notations: if $A= \big\{( x_k,y_k):k\in [\ell]\big\}\subset  {\mathbb M_{\infty}^{d}}$ then  $A^c= \left\{( y_k,x_k):k\in [\ell]\right\}\subset \Real_{+}^{d}$. 

\begin{prop}\label{Formula}Suppose that $A\subset \mathbb M_{\infty}^n $. Then for all $\bar k\in [\ell]$, we have:

$(a)$ $ \mathbb D_{\mathrm{out}}(x_{\bar k}, y_{\bar k}, T_{\mathbb Q_\infty, V,}(A))=\mathbb D_{\mathrm{in}}(-y_{\bar k}, - x_{\bar k}, T_{\mathbb Q_{-\infty}, V,}(A^c))$;

$(b)$ $  \mathbb D_{\mathrm{in}}(x_{\bar k}, y_{\bar k}, T_{\mathbb Q_\infty, V,}(A))=\mathbb D_{\mathrm{out}}(-y_{\bar k}, - x_{\bar k}, T_{\mathbb Q_{-\infty}, V,}(A^c))$;

$(c) $ $\mathbb D_{\mathrm{out}}(x_{\bar k}, y_{\bar k}, T_{\mathbb Q_{-\infty}, V,}(A))=\mathbb D_{\mathrm{in}}(-y_{\bar k}, - x_{\bar k}, T_{\mathbb Q_{+\infty}, V,}(A^c))$;

$(d)$ $ \mathbb D_{\mathrm{in}}(x_{\bar k}, y_{\bar k}, T_{\mathbb Q_{-\infty}, V,}(A))=\mathbb D_{\mathrm{out}}(-y_{\bar k}, - x_{\bar k}, T_{\mathbb Q_{+\infty}, V,}(A^c)).$

\end{prop}
{\bf Proof:} $(a)$ By definition we have, setting $s_k=-t_k$ for all $k$:

\begin{align*}\mathbb D_{\mathrm{out}}&(x_{\bar k}, y_{\bar k}, T_{\mathbb Q_\infty, V,}(A))\\& =
\max\Big\{\delta: x_{\bar k} \geq \bigvee_{k\in \ell}(t_k\un_m + x_k),  y_{\bar k}+\delta 1\!\!1_n\leq \bigvee_{k\in \ell}(t_k\un_n + y_k), \max_{k\in \ell} t_k=0\Big\}\\
&=
\max\Big\{\delta: -( -x_{\bar k} )\geq -\bigwedge_{k\in \ell}(-t_k\un_m - x_k),  -(-y_{\bar k})-(-\delta) 1\!\!1_n\leq -\bigwedge_{k\in \ell}(-t_k\un_n - y_k), \max_{k\in \ell}(- t_k)=0\Big\}\\
&=
\max\Big\{\delta: -x_{\bar k} \leq  \bigwedge_{k\in \ell}(s_k\un_m - x_k), - y_{\bar k}-\delta 1\!\!1_n\geq  \bigwedge_{k\in \ell}(s_k\un_n - y_k), \max_{k\in \ell} s_k =0\Big\}\\
&=
\mathbb D_{\mathrm{in}}(-y_{\bar k}, - x_{\bar k}, T_{\mathbb Q_{-\infty}, V,}(A^c))
  \end{align*}
$(b)$ We have
  \begin{align*}\mathbb D_{\mathrm{in}}&(x_{\bar k}, y_{\bar k}, T_{\mathbb Q_\infty, V,}(A))\\& =
\max\Big\{\delta: x_{\bar k} -\delta 1\!\!1_m\geq \bigvee_{k\in \ell}(t_k\un_m + x_k),  y_{\bar k}\leq \bigvee_{k\in \ell}(t_k\un_n + y_k), \max_{k\in \ell} t_k=0\Big\}\\
&=
\max\Big\{\delta: -( -x_{\bar k} )-(-\delta) 1\!\!1_m\geq -\bigwedge_{k\in \ell}(-t_k\un_m - x_k),  -(-y_{\bar k})\leq -\bigwedge_{k\in \ell}(-t_k\un_n - y_k), \max_{k\in \ell}(- t_k)=0\Big\}\\
&=
\max\Big\{\delta: -x_{\bar k} +\delta 1\!\!1_m\leq  \bigwedge_{k\in \ell}(s_k\un_m - x_k), - y_{\bar k}\geq  \bigwedge_{k\in \ell}(s_k\un_n - y_k), \max_{k\in \ell} s_k =0\Big\}\\
&=
\mathbb D_{\mathrm{out}}(-y_{\bar k}, - x_{\bar k}, T_{\mathbb Q_{-\infty}, V,}(A^c)).
  \end{align*}
  The proofs of $(c)$ and $(d)$ are similar. $\Box$\\

\begin{cor}Suppose that $A\subset  {\mathbb M_{\infty}^{d}} $. Then for all $\bar k\in [\ell]$, we have:

$(a)$ $\mathbb D_{\mathrm{out}}(x_{\bar k}, y_{\bar k}, T_{\mathbb Q_\infty, C,}(A))=\mathbb D_{\mathrm{in}}(-y_{\bar k}, - x_{\bar k}, T_{\mathbb Q_{-\infty}, C,}(A^c))$;

  $(b)$ $ \mathbb D_{\mathrm{in}}(x_{\bar k}, y_{\bar k}, T_{\mathbb Q_\infty, C,}(A))=\mathbb D_{\mathrm{out}}(-y_{\bar k}, - x_{\bar k}, T_{\mathbb Q_{-\infty}, C,}(A^c))$;
  
 $(c)$ $\mathbb D_{\mathrm{out}}(x_{\bar k}, y_{\bar k}, T_{\mathbb Q_{-\infty}, C,}(A))=\mathbb D_{\mathrm{in}}(-y_{\bar k}, - x_{\bar k}, T_{\mathbb Q_{+\infty}, C,}(A^c))$;
 
  $(d)$ $\mathbb D_{\mathrm{in}}(x_{\bar k}, y_{\bar k}, T_{\mathbb Q_{-\infty}, C,}(A))=\mathbb D_{\mathrm{out}}(-y_{\bar k}, - x_{\bar k}, T_{\mathbb Q_{+\infty}, C,}(A^c))$.

\end{cor}
{\bf Proof:} The proof is similar to the proof of Proposition \ref{Formula}. In particular, if $t\in \mathbb M_\infty^\ell$ then $s=-t\in \mathbb M_{-\infty}^\ell$. $\Box$\\

\begin{prop} \label{GaugeMinPolyVRS}
Let $A=\big\{( x_k,y_k):k\in [\ell]\big\}\subset \mathbb M_{- \infty}^{d}$. For all $ \bar k\in [\ell]$ such that $\beta_{\bar k,k}^c=\min\limits_{j=1,...,n}\{{y_{k,j }-y_{\bar k,j}}\}$, the input/output
translation distance functions are the following.\\ 

$(a)$ If $T=T_{\mathbb Q_{ -\infty}, C}$ then:
\begin{align}
&\mathbb D_{\mathrm{in}}\big(x_{\bar k},y_{\bar k}, T_{\mathbb Q_{ -\infty},V} (A)\big)=
\min_{i=1,\ldots,m}\max_{\substack{k\in [\ell]}}\left\{{x_{\bar
k,i} - x_{k,i}}+\min \big\{\beta^c_{\bar k,k}, 0\big\}\right\} \notag \\
&\mathbb D_{\mathrm{out}}\big(x_{\bar k},y_{\bar k}, T_{\mathbb Q_{- \infty}, V} (A) \big)=\min\Big\{\min_{i=1,\ldots,m}
\max_{\substack{k\in [\ell] \\ x_{\bar k, i}  \geq x_{k,i} }}\big\{
+x_{\bar k, ji}^{}-y_{k,j}^{} + \beta_{\bar k,k}\big\},
\max_{k\in [\ell]}\beta^c_{\bar k,k}\Big\}\notag
\end{align}

$(b)$ If
 $ T=T_{\mathbb Q_{ -\infty}, C}$ then : {\small $$  \mathbb D_{\mathrm{in}}\big (x_{\bar k},y_{\bar
k}, T_{\mathbb Q_{ -\infty},C} (A)\big)= \mathbb D_{\mathrm{out}}\big(x_{\bar k},y_{\bar k}, T_{\mathbb Q_{ -\infty}, C} (A)\big)
= \min_{i=1,...,m} \max_{k\in [\ell]}\big\{  x_{\bar k,i}^{}-x_{k,i}^{} + \beta^c_{\bar
k,k}\big\}.$$}
\end{prop}

\subsection{Discrete Production Model}

In this subsection it is assumed again that $A=\{(x_1,y_1),\ldots,(x_\ell, y_\ell)\}\subset \mathbb N^d$. This means that the components of each observed production vector are non-negative integers.  Accordingly a specific class of discrete production sets is introduced. Paralleling our earlier definitions of the Max-Plus model  proposed in \cite{abs17} lest us consider the production set defined as:
   
\begin{equation}\notag
Z_{\mathbb Q_{\infty},V} (A)  =T_{\mathbb Q_{\infty},V} (A)\cap \mathbb N^d
 \end{equation}
The graph translation homothetic model is defined as:
\begin{equation}\notag
Z_{\mathbb Q_\infty, C} (A)  =T_{\mathbb Q_\infty, C} (A) \cap \mathbb N^d
\end{equation}
It is therefore immediate to consider the lower-dequantization that is obtained taking the limit when $\alpha\longrightarrow -\infty$. It is defined as: 
\begin{equation}\notag
Z_{\mathbb Q_{-\infty}, V}(A)   =T_{\mathbb Q_{-\infty}, V}(A)  \cap \mathbb N^d
 \end{equation}
The graph translation homothetic model is defined as:

\begin{equation}\notag
Z_{\mathbb Q_{-\infty}, C}(A )   =T_{\mathbb Q_{-\infty}, C}(A)  \cap \mathbb N^d
 \end{equation}

\bigskip 
\begin{center}{\scriptsize  
 \unitlength 0.4mm 
\linethickness{0.4pt}
\ifx\plotpoint\undefined\newsavebox{\plotpoint}\fi 
\begin{picture}(269,120.5)(0,0)
\put(4,11.5){\vector(0,1){100.75}}
\put(159.75,11.75){\vector(0,1){100.75}}
\put(4,11.25){\vector(1,0){102.5}}
\put(159.75,11.5){\vector(1,0){102.5}}
\put(113.25,11){\makebox(0,0)[cc]{$x$}}
\put(269,11.25){\makebox(0,0)[cc]{$x$}}
\put(36.75,66.25){\makebox(0,0)[cc]{$z_1$}}
\put(192.5,66.5){\makebox(0,0)[cc]{$z_1$}}
\put(100.75,118){\makebox(0,0)[cc]{$z_2$}}
\put(255.75,115.25){\makebox(0,0)[cc]{$z_2$}}
\put(70.25,52.75){\makebox(0,0)[cc]{$z_3$}}
\put(226,53){\makebox(0,0)[cc]{$z_3$}}
\put(25.5,29.5){\makebox(0,0)[cc]{$z_5$}}
\put(181.25,29.75){\makebox(0,0)[cc]{$z_5$}}
\put(54.75,102.5){\makebox(0,0)[cc]{$z_4$}}
\put(210.5,102.75){\makebox(0,0)[cc]{$z_4$}}
\put(3.75,120.25){\makebox(0,0)[cc]{$y$}}
\put(159.5,120.5){\makebox(0,0)[cc]{$y$}}
\put(0,7.75){\makebox(0,0)[cc]{$0$}}
\put(155.75,8){\makebox(0,0)[cc]{$0$}}
\put(69.5,0){\makebox(0,0)[cc]
{{{\bf Figure \ref{Dequantized}.5:} Max-Plus Discrete Production Set.\ \ \ }}}
\put(225.25,.25){\makebox(0,0)[cc]
{{{\bf Figure \ref{Dequantized}.6:} Min-Plus Discrete Production Set.}}}
\put(3.75,11.75){\vector(1,1){9.25}}
\put(159.5,12){\vector(1,1){9.25}}
\put(82.25,67.5){\makebox(0,0)[cc]{$Z_{\mathbb Q_\infty,V}$}}
\put(238,67.75){\makebox(0,0)[cc]{$Z_{\mathbb Q_{-\infty},V}$}}
\put(17,80.75){\circle*{1.5}}
\put(172.75,81){\circle*{1.5}}
\put(16.75,49.5){\circle*{1.5}}
\put(172.5,49.75){\circle*{1.5}}
\put(17.25,111.25){\circle*{1.5}}
\put(173,111.5){\circle*{1.5}}
\put(30.25,80.75){\circle*{1.5}}
\put(186,81){\circle*{1.5}}
\put(30,49.5){\circle*{1.5}}
\put(185.75,49.75){\circle*{1.5}}
\put(30.5,111.25){\circle*{1.5}}
\put(186.25,111.5){\circle*{1.5}}
\put(57.25,80.25){\circle*{1.5}}
\put(213,80.5){\circle*{1.5}}
\put(57,49){\circle*{1.5}}
\put(212.75,49.25){\circle*{1.5}}
\put(57.5,110.75){\circle*{1.5}}
\put(213.25,111){\circle*{1.5}}
\put(83.75,80.75){\circle*{1.5}}
\put(239.5,81){\circle*{1.5}}
\put(83.5,49.5){\circle*{1.5}}
\put(239.25,49.75){\circle*{1.5}}
\put(84,111.25){\circle*{1.5}}
\put(239.75,111.5){\circle*{1.5}}
\put(16.75,67){\circle*{1.5}}
\put(172.5,67.25){\circle*{1.5}}
\put(16.5,35.75){\circle*{1.5}}
\put(172.25,36){\circle*{1.5}}
\put(17,24){\circle*{1.5}}
\put(172.75,24.25){\circle*{1.5}}
\put(17,97.5){\circle*{1.5}}
\put(172.75,97.75){\circle*{1.5}}
\put(30,67){\circle*{1.5}}
\put(185.75,67.25){\circle*{1.5}}
\put(29.75,35.75){\circle*{1.5}}
\put(185.5,36){\circle*{1.5}}
\put(30.25,24){\circle*{1.5}}
\put(186,24.25){\circle*{1.5}}
\put(30.25,97.5){\circle*{1.5}}
\put(186,97.75){\circle*{1.5}}
\put(57,66.5){\circle*{1.5}}
\put(212.75,66.75){\circle*{1.5}}
\put(56.75,35.25){\circle*{1.5}}
\put(212.5,35.5){\circle*{1.5}}
\put(57.25,23.5){\circle*{1.5}}
\put(213,23.75){\circle*{1.5}}
\put(57.25,97){\circle*{1.5}}
\put(213,97.25){\circle*{1.5}}
\put(83.5,67){\circle*{1.5}}
\put(239.25,67.25){\circle*{1.5}}
\put(83.25,35.75){\circle*{1.5}}
\put(239,36){\circle*{1.5}}
\put(83.75,24){\circle*{1.5}}
\put(239.5,24.25){\circle*{1.5}}
\put(83.75,97.5){\circle*{1.5}}
\put(239.5,97.75){\circle*{1.5}}
\put(30.25,80.5){\circle*{1.5}}
\put(186,80.75){\circle*{1.5}}
\put(30,49.25){\circle*{1.5}}
\put(185.75,49.5){\circle*{1.5}}
\put(30.5,111){\circle*{1.5}}
\put(186.25,111.25){\circle*{1.5}}
\put(43.5,80.5){\circle*{1.5}}
\put(199.25,80.75){\circle*{1.5}}
\put(43.25,49.25){\circle*{1.5}}
\put(199,49.5){\circle*{1.5}}
\put(43.75,111){\circle*{1.5}}
\put(199.5,111.25){\circle*{1.5}}
\put(70.5,80){\circle*{1.5}}
\put(226.25,80.25){\circle*{1.5}}
\put(70.25,48.75){\circle*{1.5}}
\put(226,49){\circle*{1.5}}
\put(70.75,110.5){\circle*{1.5}}
\put(226.5,110.75){\circle*{1.5}}
\put(97,80.5){\circle*{1.5}}
\put(252.75,80.75){\circle*{1.5}}
\put(96.75,49.25){\circle*{1.5}}
\put(252.5,49.5){\circle*{1.5}}
\put(97.25,111){\circle*{1.5}}
\put(253,111.25){\circle*{1.5}}
\put(30,66.75){\circle*{1.5}}
\put(185.75,67){\circle*{1.5}}
\put(29.75,35.5){\circle*{1.5}}
\put(185.5,35.75){\circle*{1.5}}
\put(30.25,23.75){\circle*{1.5}}
\put(186,24){\circle*{1.5}}
\put(30.25,97.25){\circle*{1.5}}
\put(186,97.5){\circle*{1.5}}
\put(43.25,66.75){\circle*{1.5}}
\put(199,67){\circle*{1.5}}
\put(43,35.5){\circle*{1.5}}
\put(198.75,35.75){\circle*{1.5}}
\put(43.5,23.75){\circle*{1.5}}
\put(199.25,24){\circle*{1.5}}
\put(43.5,97.25){\circle*{1.5}}
\put(199.25,97.5){\circle*{1.5}}
\put(70.25,66.25){\circle*{1.5}}
\put(226,66.5){\circle*{1.5}}
\put(70,35){\circle*{1.5}}
\put(225.75,35.25){\circle*{1.5}}
\put(70.5,23.25){\circle*{1.5}}
\put(226.25,23.5){\circle*{1.5}}
\put(70.5,96.75){\circle*{1.5}}
\put(226.25,97){\circle*{1.5}}
\put(96.75,66.75){\circle*{1.5}}
\put(252.5,67){\circle*{1.5}}
\put(96.5,35.5){\circle*{1.5}}
\put(252.25,35.75){\circle*{1.5}}
\put(97,23.75){\circle*{1.5}}
\put(252.75,24){\circle*{1.5}}
\put(97,97.25){\circle*{1.5}}
\put(252.75,97.5){\circle*{1.5}}
\put(30,17.5){\line(0,1){33}}
\multiput(30,50.5)(.083870968,.103225806){155}{\line(0,1){.103225806}}
\put(43,66.5){\line(0,1){14.75}}
\multiput(43,81.25)(.084319527,.093195266){169}{\line(0,1){.093195266}}
\multiput(97.25,111.5)(-.085403727,-.083850932){161}{\line(-1,0){.085403727}}
\multiput(83.25,97.75)(-8.66667,-.08333){3}{\line(-1,0){8.66667}}
\put(97.5,111.5){\line(1,0){7.25}}
\put(185.75,36.25){\line(0,-1){19}}
\multiput(186,36.75)(.083881579,.085526316){152}{\line(0,1){.085526316}}
\put(198.75,49.75){\line(0,1){17.75}}
\multiput(198.75,67.5)(.088607595,.083860759){158}{\line(1,0){.088607595}}
\put(212.75,80.75){\line(0,1){16.5}}
\multiput(212.75,97.25)(.083841463,.086890244){164}{\line(0,1){.086890244}}
\put(226.5,111.5){\line(1,0){32}}
\end{picture}

}
 \end{center}
 \bigskip
 
In Propositions  \ref{GaugeMaxPolyVRS} and  \ref{GaugeMinPolyVRS}, the computation of input and output oriented translation functions only involves addition, difference and the computation of minima and maxima. Since the components of the production vector are assumed to be integer valued, the translation distance functions are also integer valued. 
 
 \begin{prop} Suppose that  {$A=\{(x_k,y_k):k\in [\ell]\}\subset  \mathbb N^d$}. Let us consider the collection   $\mathcal Z(A)$ defined as:
 $$\mathcal T (A)=\{T_{\mathbb Q_{\infty},V} (A), T_{\mathbb Q_{\infty},C} (A),T_{\mathbb Q_{-\infty},V} (A), T_{\mathbb Q_{-\infty},C} (A) \}.$$
 Then, for all $T\in \mathcal T(A)$, and all $\bar k\in [\ell]$ we have:
 $$\mathbb D_{\mathrm{in}}(x_{\bar k},y_{\bar k}, T)\in \mathbb N\quad \text{ and }\quad \mathbb D_{\mathrm{out}}(x,y, T)\in \mathbb N. $$

 \end{prop}
 {\bf Proof:} Since  $A\subset \mathbb N^d  $, it follows that $\beta_{\bar k,k}=\min_{i=1,\ldots,m}\{{x_i^{\bar k}}-x_i^k \} \in  \mathbb Z$. Moreover for any $i,j,k$ $x_{k,i}\in \mathbb N$ and $y_{j,k}\in \mathbb N$. Therefore from Propositions  \ref{GaugeMaxPolyVRS} and  \ref{GaugeMinPolyVRS} and it follows that
 $$\mathbb D_{\mathrm{in}}\big(x_{\bar k},y_{\bar k}, T_{\mathbb Q_{ \infty}, V}(A )  \big)\in \mathbb Z\quad \text{ and }\quad \mathbb D_{\mathrm{out}}\big (x_{\bar k},y_{\bar k}, T_{\mathbb Q_{ \infty}, V}(A )\big)\in \mathbb Z. $$
 However $\mathbb D_{\mathrm{in}}\big (x_{\bar k},y_{\bar k}, T_{\mathbb Q_{ \infty}, V}(A )  \big)\geq 0$ and $\mathbb D_{\mathrm{out}}\big(x_{\bar k},y_{\bar k}, T_{\mathbb Q_{ \infty}, V}(A ) \big )  {\geq 0}$. Hence $\mathbb D_{\mathrm{in}}\big (x_{\bar k},y_{\bar k}, T_{\mathbb Q_{ \infty}, V}(A ) \big )\in \mathbb N$ and $\mathbb D_{\mathrm{out}}\big (x_{\bar k},y_{\bar k}, T_{\mathbb Q_{ \infty}, V}(A )\big )\in \mathbb N. $ The proof is similar for all $T\in \mathcal T (A).$ $\Box$\\

 \begin{prop} For all finite subset $A=\{(x_k,y_k):k\in [\ell]\}\subset  \mathbb N^d$, the following relations arise: 
 
 $(a)$ $\mathbb D_{\mathrm{in}}(x_{\bar k},y_{\bar k}, T_{\mathbb Q_\infty, V} )=\mathbb D_{\mathrm{in}}(x_{\bar k},y_{\bar k}, Z_{\mathbb Q_\infty, V} )$;
 
  $(b)$ $\mathbb D_{\mathrm{in}}(x_{\bar k},y_{\bar k}, T_{\mathbb Q_\infty,C} )=\mathbb D_{\mathrm{in}}(x_{\bar k},y_{\bar k}, Z_{\mathbb Q_\infty, C} )$;
  
   $(c)$ $\mathbb D_{\mathrm{out}}(x_{\bar k},y_{\bar k}, T_{\mathbb Q_\infty, V} )=\mathbb D_{\mathrm{out}}(x_{\bar k},y_{\bar k}, Z_{\mathbb Q_\infty, V} )$;
   
  $(d)$ $\mathbb D_{\mathrm{out}}(x_{\bar k},y_{\bar k}, T_{\mathbb Q_\infty,C} )=\mathbb D_{\mathrm{out}}(x_{\bar k},y_{\bar k}, Z_{\mathbb Q_\infty, C} )$.

 \end{prop}
 {\bf Proof:}  $(a)$ By hypothesis $Z_{\mathbb Q_\infty, V}\subset  T_{\mathbb Q_\infty, V}$. Therefore:
$$
\mathbb D_{\mathrm{in}}\big(x_{\bar k},y_{\bar k}, T_{\mathbb Q_\infty, V}  \big) \geq \mathbb D_{\mathrm{in}}\big(x_{\bar k},y_{\bar k}, Z_{\mathbb Q_\infty, V}  \big)  
$$
  However, we have proved that  for all $T\in \mathcal T $, and all $\bar k\in [\ell]$ we have:
$$
\mathbb D_{\mathrm{in}}(x_{\bar k},y_{\bar k}, T)\in \mathbb N $$
 Therefore, it follows that:
 $$\Big(x_{\bar k}-\mathbb D_{\mathrm{in}}(x_{\bar k},y_{\bar k}, T_{\mathbb Q_\infty, V})\un_m,y_{\bar k} \Big)\in \mathbb N^d 
$$
Consequently,
$$
\Big(x_{\bar k}-\mathbb D_{\mathrm{in}}\big(x_{\bar k},y_{\bar k}, T_{\mathbb Q_\infty, V} \big)\un_m,y_{\bar k} \Big)\in T_{\mathbb Q_\infty, V}\cap \mathbb N^d=Z_{\mathbb Q_\infty, V}
$$
Then,
$$
\mathbb D_{\mathrm{in}}\big (x_{\bar k},y_{\bar k}, T_{\mathbb Q_\infty, V}  \big )\leq \mathbb D_{\mathrm{in}}\big(x_{\bar k},y_{\bar k}, Z_{\mathbb Q_\infty, V} \big) 
$$
It follows that:    $\mathbb D_{\mathrm{in}}\big (x_{\bar k},y_{\bar k}, T_{\mathbb Q_\infty, V}  \big )= \mathbb D_{\mathrm{in}}\big (x_{\bar k},y_{\bar k}, Z_{\mathbb Q_\infty, V}  \big ). $
  The proof is similar for  $(b)$, $(c)$ and $(d)$. Consequently the result holds for all $T\in \mathcal T  $. $\Box$\\

\begin{center}{\scriptsize 
\unitlength 0.4mm 
\linethickness{0.4pt}
\ifx\plotpoint\undefined\newsavebox{\plotpoint}\fi 
\begin{picture}(269,120.5)(0,0)
\put(4,11.5){\vector(0,1){100.75}}
\put(159.75,11.75){\vector(0,1){100.75}}
\put(4,11.25){\vector(1,0){102.5}}
\put(159.75,11.5){\vector(1,0){102.5}}
\put(113.25,11){\makebox(0,0)[cc]{$x$}}
\put(269,11.25){\makebox(0,0)[cc]{$x$}}
\put(3.75,120.25){\makebox(0,0)[cc]{$y$}}
\put(159.5,120.5){\makebox(0,0)[cc]{$y$}}
\put(0,7.75){\makebox(0,0)[cc]{$0$}}
\put(155.75,8){\makebox(0,0)[cc]{$0$}}
\put(50.5,0){\makebox(0,0)[cc]
{{{\bf Figure \ref{Dequantized}.7:} Max-Plus Discrete Input distance Function. \ \ \ \ \ \ }}}
\put(235.25,.25){\makebox(0,0)[cc]
{{{\bf Figure \ref{Dequantized}.8:} Max-Plus Discrete Output Distance Function.}}}
\put(3.75,11.75){\vector(1,1){9.25}}
\put(159.5,12){\vector(1,1){9.25}}
\put(17,80.75){\circle*{1.5}}
\put(172.75,81){\circle*{1.5}}
\put(16.75,49.5){\circle*{1.5}}
\put(172.5,49.75){\circle*{1.5}}
\put(17.25,111.25){\circle*{1.5}}
\put(173,111.5){\circle*{1.5}}
\put(30.25,80.75){\circle*{1.5}}
\put(186,81){\circle*{1.5}}
\put(30,49.5){\circle*{1.5}}
\put(185.75,49.75){\circle*{1.5}}
\put(30.5,111.25){\circle*{1.5}}
\put(186.25,111.5){\circle*{1.5}}
\put(213,80.5){\circle*{1.5}}
\put(57,49){\circle*{1.5}}
\put(212.75,49.25){\circle*{1.5}}
\put(57.5,110.75){\circle*{1.5}}
\put(213.25,111){\circle*{1.5}}
\put(83.75,80.75){\circle*{1.5}}
\put(239.5,81){\circle*{1.5}}
\put(83.5,49.5){\circle*{1.5}}
\put(239.25,49.75){\circle*{1.5}}
\put(84,111.25){\circle*{1.5}}
\put(239.75,111.5){\circle*{1.5}}
\put(16.75,67){\circle*{1.5}}
\put(172.5,67.25){\circle*{1.5}}
\put(16.5,35.75){\circle*{1.5}}
\put(172.25,36){\circle*{1.5}}
\put(17,24){\circle*{1.5}}
\put(172.75,24.25){\circle*{1.5}}
\put(17,97.5){\circle*{1.5}}
\put(172.75,97.75){\circle*{1.5}}
\put(30,67){\circle*{1.5}}
\put(185.75,67.25){\circle*{1.5}}
\put(29.75,35.75){\circle*{1.5}}
\put(185.5,36){\circle*{1.5}}
\put(30.25,24){\circle*{1.5}}
\put(186,24.25){\circle*{1.5}}
\put(30.25,97.5){\circle*{1.5}}
\put(186,97.75){\circle*{1.5}}
\put(57,66.5){\circle*{1.5}}
\put(212.75,66.75){\circle*{1.5}}
\put(56.75,35.25){\circle*{1.5}}
\put(212.5,35.5){\circle*{1.5}}
\put(57.25,23.5){\circle*{1.5}}
\put(213,23.75){\circle*{1.5}}
\put(57.25,97){\circle*{1.5}}
\put(56.75,80.25){\circle*{1.5}}
\put(213,97.25){\circle*{1.5}}
\put(83.5,67){\circle*{1.5}}
\put(239.25,67.25){\circle*{1.5}}
\put(83.25,35.75){\circle*{1.5}}
\put(239,36){\circle*{1.5}}
\put(83.75,24){\circle*{1.5}}
\put(239.5,24.25){\circle*{1.5}}
\put(83.75,97.5){\circle*{1.5}}
\put(239.5,97.75){\circle*{1.5}}
\put(30.25,80.5){\circle*{1.5}}
\put(186,80.75){\circle*{1.5}}
\put(30,49.25){\circle*{1.5}}
\put(185.75,49.5){\circle*{1.5}}
\put(30.5,111){\circle*{1.5}}
\put(186.25,111.25){\circle*{1.5}}
\put(43.5,80.5){\circle*{1.5}}
\put(199.25,81.5){\circle*{1.5}}
\put(43.25,49.25){\circle*{1.5}}
\put(199,49.5){\circle*{1.5}}
\put(43.75,111){\circle*{1.5}}
\put(199.5,111.25){\circle*{1.5}}
\put(70.5,80){\circle*{1.5}}
\put(226.25,80.25){\circle*{1.5}}
\put(70.25,48.75){\circle*{1.5}}
\put(226,49){\circle*{1.5}}
\put(70.75,110.5){\circle*{1.5}}
\put(226.5,110.75){\circle*{1.5}}
\put(97,80.5){\circle*{1.5}}
\put(252.75,80.75){\circle*{1.5}}
\put(96.75,49.25){\circle*{1.5}}
\put(252.5,49.5){\circle*{1.5}}
\put(97.25,111){\circle*{1.5}}
\put(253,111.25){\circle*{1.5}}
\put(30,66.75){\circle*{1.5}}
\put(185.75,67){\circle*{1.5}}
\put(29.75,35.5){\circle*{1.5}}
\put(185.5,35.75){\circle*{1.5}}
\put(30.25,23.75){\circle*{1.5}}
\put(186,24){\circle*{1.5}}
\put(30.25,97.25){\circle*{1.5}}
\put(186,97.5){\circle*{1.5}}
\put(44,66.5){\circle*{1.5}}
\put(199,67){\circle*{1.5}}
\put(43,35.5){\circle*{1.5}}
\put(198.75,35.75){\circle*{1.5}}
\put(43.5,23.75){\circle*{1.5}}
\put(199.25,24){\circle*{1.5}}
\put(43.5,97.25){\circle*{1.5}}
\put(199.25,97.5){\circle*{1.5}}
\put(70.25,66.25){\circle*{1.5}}
\put(226,66.5){\circle*{1.5}}
\put(70,35){\circle*{1.5}}
\put(225.75,35.25){\circle*{1.5}}
\put(70.5,23.25){\circle*{1.5}}
\put(226.25,23.5){\circle*{1.5}}
\put(70.5,96.75){\circle*{1.5}}
\put(226.25,97){\circle*{1.5}}
\put(96.75,66.75){\circle*{1.5}}
\put(252.5,67){\circle*{1.5}}
\put(96.5,35.5){\circle*{1.5}}
\put(252.25,35.75){\circle*{1.5}}
\put(97,23.75){\circle*{1.5}}
\put(252.75,24){\circle*{1.5}}
\put(97,97.25){\circle*{1.5}}
\put(252.75,97.5){\circle*{1.5}}
\put(17.25,111.75){\line(0,-1){44}}
\put(17.25,67.75){\line(1,0){12.75}}
\put(30,67.75){\line(0,-1){18}}
\put(30,49.75){\line(1,0){14.25}}
\put(43,48.25){\line(0,-1){13.25}}
\put(43,35){\line(1,0){13.5}}
\put(56.5,35){\line(0,-1){11.75}}
\put(56.5,23.25){\line(1,0){39.75}}
\put(29.5,49){\vector(-3,-4){.176}}\multiput(70.25,96.5)(-.0841942149,-.0981404959){484}{\line(0,-1){.0981404959}}
\put(73,89.5){\makebox(0,0)[cc]{$x_{\bar k}$}}
\multiput(160,97.25)(8.833333,.083333){9}{\line(1,0){8.833333}}
\multiput(239.5,98)(-.08333,-28.25){3}{\line(0,-1){28.25}}
\put(213.25,97.75){\vector(1,1){.176}}\multiput(186,67.75)(.0841049383,.0925925926){324}{\line(0,1){.0925925926}}
\put(186.25,59){\makebox(0,0)[cc]{$y_{\bar k}$}}
\end{picture}}
\end{center}
\bigskip

The proof of the next statement is similar.
 
\begin{prop} Suppose that $A=\{(x_k,y_k):k\in [\ell]\}\subset  \mathbb N^d$.  The following relations arise:
 
 $(a)$ $\mathbb D_{\mathrm{in}}(x_{\bar k},y_{\bar k}, T_{\mathbb Q_{-\infty}, V} )=\mathbb D_{\mathrm{in}}(x_{\bar k},y_{\bar k}, Z_{\mathbb Q_{-\infty}, V} )$;
 
  $(b)$ $\mathbb D_{\mathrm{in}}(x_{\bar k},y_{\bar k}, T_{\mathbb Q_{-\infty},C} )=\mathbb D_{\mathrm{in}}(x_{\bar k},y_{\bar k}, Z_{\mathbb Q_{-\infty}, C} )$;
  
   $(c)$ $\mathbb D_{\mathrm{out}}(x_{\bar k},y_{\bar k}, T_{\mathbb Q_{-\infty}, V} )=\mathbb D_{\mathrm{out}}(x_{\bar k},y_{\bar k}, Z_{\mathbb Q_{-\infty}, V} )$;
   
  $(d)$ $\mathbb D_{\mathrm{out}}(x_{\bar k},y_{\bar k}, T_{\mathbb Q_{-\infty},C} )=\mathbb D_{\mathrm{out}}(x_{\bar k},y_{\bar k}, Z_{\mathbb Q_{-\infty}, C} )$.
\end{prop}

\bigskip

\begin{center}
{\scriptsize 

\unitlength 0.4mm 
\linethickness{0.4pt}
\ifx\plotpoint\undefined\newsavebox{\plotpoint}\fi 
\begin{picture}(255.25,120.5)(0,0)
\put(146,11.5){\vector(0,1){100.75}}
\put(4,11.75){\vector(0,1){100.75}}
\put(146,11.25){\vector(1,0){102.5}}
\put(4,11.5){\vector(1,0){102.5}}
\put(255.25,11){\makebox(0,0)[cc]{$x$}}
\put(113.25,11.25){\makebox(0,0)[cc]{$x$}}
\put(145.75,120.25){\makebox(0,0)[cc]{$y$}}
\put(3.75,120.5){\makebox(0,0)[cc]{$y$}}
\put(142,7.75){\makebox(0,0)[cc]{$0$}}
\put(0,8){\makebox(0,0)[cc]{$0$}}
\put(50.5,.25){\makebox(0,0)[cc]
{{{\bf Figure \ref{Dequantized}.9:} Min-Plus  Discrete Input Distance Function.\ \ \ \ \ \ \ }}}
\put(236.529,.25){\makebox(0,0)[cc]
{{{\bf Figure \ref{Dequantized}.10:} Min-Plus  Discrete Output Distance Function.}}}
\put(145.75,11.75){\vector(1,1){9.25}}
\put(3.75,12){\vector(1,1){9.25}}
\put(159,80.75){\circle*{1.5}}
\put(17,81){\circle*{1.5}}
\put(158.75,49.5){\circle*{1.5}}
\put(16.75,49.75){\circle*{1.5}}
\put(159.25,111.25){\circle*{1.5}}
\put(17.25,111.5){\circle*{1.5}}
\put(172.25,80.75){\circle*{1.5}}
\put(30.25,81){\circle*{1.5}}
\put(172,49.5){\circle*{1.5}}
\put(30,49.75){\circle*{1.5}}
\put(172.5,111.25){\circle*{1.5}}
\put(30.5,111.5){\circle*{1.5}}
\put(57.25,80.5){\circle*{1.5}}
\put(199,49){\circle*{1.5}}
\put(57,49.25){\circle*{1.5}}
\put(199.5,110.75){\circle*{1.5}}
\put(57.5,111){\circle*{1.5}}
\put(225.75,80.75){\circle*{1.5}}
\put(83.75,81){\circle*{1.5}}
\put(225.5,49.5){\circle*{1.5}}
\put(83.5,49.75){\circle*{1.5}}
\put(226,111.25){\circle*{1.5}}
\put(84,111.5){\circle*{1.5}}
\put(158.75,67){\circle*{1.5}}
\put(16.75,67.25){\circle*{1.5}}
\put(158.5,35.75){\circle*{1.5}}
\put(16.5,36){\circle*{1.5}}
\put(159,24){\circle*{1.5}}
\put(17,24.25){\circle*{1.5}}
\put(159,97.5){\circle*{1.5}}
\put(17,97.75){\circle*{1.5}}
\put(172,67){\circle*{1.5}}
\put(30,67.25){\circle*{1.5}}
\put(172.25,35){\circle*{1.5}}
\put(29.75,36){\circle*{1.5}}
\put(172.25,24){\circle*{1.5}}
\put(30.25,24.25){\circle*{1.5}}
\put(172.25,97.5){\circle*{1.5}}
\put(30.25,97.75){\circle*{1.5}}
\put(199,66.5){\circle*{1.5}}
\put(57,66.75){\circle*{1.5}}
\put(198.75,35.25){\circle*{1.5}}
\put(56.75,35.5){\circle*{1.5}}
\put(199.25,23.5){\circle*{1.5}}
\put(57.25,23.75){\circle*{1.5}}
\put(199.25,97){\circle*{1.5}}
\put(198.75,80.25){\circle*{1.5}}
\put(57.25,97.25){\circle*{1.5}}
\put(225.5,67){\circle*{1.5}}
\put(83.5,66.358){\circle*{1.5}}
\put(225.25,35.75){\circle*{1.5}}
\put(83.25,36){\circle*{1.5}}
\put(225.75,24){\circle*{1.5}}
\put(83.75,24.25){\circle*{1.5}}
\put(225.75,97.5){\circle*{1.5}}
\put(83.75,97.75){\circle*{1.5}}
\put(172.25,80.5){\circle*{1.5}}
\put(30.25,80.75){\circle*{1.5}}
\put(172,49.25){\circle*{1.5}}
\put(30,49.5){\circle*{1.5}}
\put(172.5,111){\circle*{1.5}}
\put(30.5,111.25){\circle*{1.5}}
\put(185.5,80.5){\circle*{1.5}}
\put(43.5,81.5){\circle*{1.5}}
\put(185.25,49.25){\circle*{1.5}}
\put(43.25,49.5){\circle*{1.5}}
\put(185.75,111){\circle*{1.5}}
\put(43.75,111.25){\circle*{1.5}}
\put(212.5,80){\circle*{1.5}}
\put(70.5,80.25){\circle*{1.5}}
\put(212.25,48.75){\circle*{1.5}}
\put(70.25,49){\circle*{1.5}}
\put(212.75,110.5){\circle*{1.5}}
\put(70.75,110.75){\circle*{1.5}}
\put(239,80.5){\circle*{1.5}}
\put(97,80.75){\circle*{1.5}}
\put(238.75,49.25){\circle*{1.5}}
\put(96.75,49.5){\circle*{1.5}}
\put(239.25,111){\circle*{1.5}}
\put(97.25,111.25){\circle*{1.5}}
\put(172,66.75){\circle*{1.5}}
\put(30,67){\circle*{1.5}}
\put(172.25,34.75){\circle*{1.5}}
\put(29.75,35.75){\circle*{1.5}}
\put(172.25,23.75){\circle*{1.5}}
\put(30.25,24){\circle*{1.5}}
\put(172.25,97.25){\circle*{1.5}}
\put(30.25,97.5){\circle*{1.5}}
\put(186,66.5){\circle*{1.5}}
\put(43.25,67){\circle*{1.5}}
\put(185,35.5){\circle*{1.5}}
\put(43,35.75){\circle*{1.5}}
\put(185.5,23.75){\circle*{1.5}}
\put(43.5,24){\circle*{1.5}}
\put(185.5,97.25){\circle*{1.5}}
\put(43.5,97.5){\circle*{1.5}}
\put(212.25,66.25){\circle*{1.5}}
\put(70.25,66.5){\circle*{1.5}}
\put(212,35){\circle*{1.5}}
\put(70,35.25){\circle*{1.5}}
\put(212.5,23.25){\circle*{1.5}}
\put(70.5,23.5){\circle*{1.5}}
\put(212.5,96.75){\circle*{1.5}}
\put(70.5,97){\circle*{1.5}}
\put(238.75,66.75){\circle*{1.5}}
\put(96.75,67){\circle*{1.5}}
\put(238.5,35.5){\circle*{1.5}}
\put(96.5,35.75){\circle*{1.5}}
\put(239,23.75){\circle*{1.5}}
\put(97,24){\circle*{1.5}}
\put(239,97.25){\circle*{1.5}}
\put(97,97.5){\circle*{1.5}}
\put(212.75,81.75){\vector(3,4){.176}}\multiput(172,34.25)(.0841942149,.0981404959){484}{\line(0,1){.0981404959}}
\put(169.25,41.25){\makebox(0,0)[cc]{$x_{\bar k}$}}
\multiput(109.321,35.824)(-8.833333,-.083333){9}{\line(-1,0){8.833333}}
\multiput(29.821,35.074)(.08333,28.25){3}{\line(0,1){28.25}}
\put(83.071,74.074){\makebox(0,0)[]{$y_{\bar k}$}}
\put(186.75,112.25){\line(0,1){0}}
\multiput(159,111.75)(8.91667,.08333){3}{\line(1,0){8.91667}}
\put(185.75,112){\line(0,-1){14.25}}
\put(185.75,97.75){\line(1,0){13.25}}
\put(199,97.75){\line(0,-1){16.5}}
\put(199,81.25){\line(1,0){13.25}}
\put(212.25,81.25){\line(0,-1){33}}
\put(212,49){\line(1,0){14}}
\put(226,49){\line(0,-1){32.75}}
\put(226,16.25){\line(0,1){0}}
\put(57.231,35.32){\vector(-1,-1){.176}}\multiput(83.394,64.753)(-.084123976,-.094639473){311}{\line(0,-1){.094639473}}
\end{picture}
}
\end{center}

  \subsection{Free Disposal Hull}

 Suppose that  $A=\{(x_k,y_k):k\in [\ell]\}\subset  \mathbb N^d$. By definition  $T_{F }(A)=(A+K)\cap \Real_{+}^{d}.$ Paralleling our earlier definitions, we can define the set:
 \begin{equation}\notag
 Z_{F }(A)=(A+K)\cap \Real_{+}^{d} 
= T_{F }(A)=(A+K)\cap \mathbb N^{d}
\end{equation}
In the following it is shown that the directional distance function, when it is either input or output oriented, has integer values. It follows that the benchmarking points also have integers components.

\begin{prop}Suppose that $A=\{(x_k,y_k):k\in [\ell]\}\subset  \mathbb N^d$. Then for all $\bar k\in [\ell]$ we have;
$$\mathbb D_{\mathrm{in}}\big(x_{\bar k}, y_{\bar k}, T_F(A)\big)\in \mathbb N\quad \text{ and }\quad \mathbb D_{\mathrm{out}}\big(x_{\bar k}, y_{\bar k}, T_F(A)\big) \in \mathbb N $$   
\end{prop}
{\bf Proof:} Suppose that $\delta=\mathbb D_{\mathrm{in}}\big(x_{\bar k}, y_{\bar k}, T_F(A)\big)$. Then there is some $i\in [m]$ such that
$$x_{\bar k,i}+\delta =x_{  k_0,i}$$
for some $k_0\in [\ell]$. However, by hypothesis $z_{\bar k},z_{k_0}\in A\subset \mathbb N^d.$ Therefore:
$$\delta =x_{  k_0,i}-x_{\bar k,i}\in  \mathbb N$$
which ends the first par of the proof. The proof of the second part is similar. $\Box$\\

 \begin{prop} For all finite subset $A=\{(x_k,y_k):k\in [\ell]\}\subset  \mathbb N^d$,  we have the relations: 
 
 $(a)$ $\mathbb D_{\mathrm{in}}(x_{\bar k},y_{\bar k}, T_{F} )=\mathbb D_{\mathrm{in}}(x_{\bar k},y_{\bar k}, Z_{F} )$;

   $(b)$ $\mathbb D_{\mathrm{out}}(x_{\bar k},y_{\bar k}, T_{F} )=\mathbb D_{\mathrm{out}}(x_{\bar k},y_{\bar k}, Z_{F} )$.

 \end{prop}
 {\bf Proof:}  $(a)$ By hypothesis $Z_{F}\subset  T_{F}$. Therefore:
$$
\mathbb D_{\mathrm{in}}\big(x_{\bar k},y_{\bar k}, T_{F}  \big) \geq \mathbb D_{\mathrm{in}}\big(x_{\bar k},y_{\bar k}, Z_{F}  \big)  
$$
However it has been shown that for all $\bar k\in [\ell]$ we have:
$$
\mathbb D_{\mathrm{in}}(x_{\bar k},y_{\bar k}, F)\in \mathbb N 
$$
Therefore, it follows that:
$$
\Big(x_{\bar k}-\mathbb D_{\mathrm{in}}(x_{\bar k},y_{\bar k}, T_{F})\un_m,y_{\bar k} \Big)\in \mathbb N^d 
$$
Then,
$$
\Big(x_{\bar k}-\mathbb D_{\mathrm{in}}\big(x_{\bar k},y_{\bar k}, T_{F} \big)\un_m,y_{\bar k} \Big)\in T_{F}\cap \mathbb N^d=Z_{F} 
$$
Consequently,   
$$
\mathbb D_{\mathrm{in}}\big (x_{\bar k},y_{\bar k}, T_{F}  \big )\leq \mathbb D_{\mathrm{in}}\big(x_{\bar k},y_{\bar k}, Z_{F} \big) 
$$
It follows that:    $\mathbb D_{\mathrm{in}}\big (x_{\bar k},y_{\bar k}, T_{F}  \big )= \mathbb D_{\mathrm{in}}\big (x_{\bar k},y_{\bar k}, Z_{F}  \big ). $
  The proof is similar for  $(b)$.  $\Box$\\

 \bigskip
 \begin{center}
 {\scriptsize 
\unitlength 0.4mm 
\linethickness{0.4pt}
\ifx\plotpoint\undefined\newsavebox{\plotpoint}\fi 
\begin{picture}(255.25,120.25)(0,0)
\put(146,11.25){\vector(0,1){100.75}}
\put(4,11.5){\vector(0,1){100.75}}
\put(146,11){\vector(1,0){102.5}}
\put(4,11.25){\vector(1,0){102.5}}
\put(255.25,10.75){\makebox(0,0)[cc]{$x$}}
\put(113.25,11){\makebox(0,0)[cc]{$x$}}
\put(145.75,120){\makebox(0,0)[cc]{$y$}}
\put(3.75,120.25){\makebox(0,0)[cc]{$y$}}
\put(142,7.5){\makebox(0,0)[cc]{$0$}}
\put(0,7.75){\makebox(0,0)[cc]{$0$}}
\put(50.5,0){\makebox(0,0)[cc]
{{{\bf Figure \ref{Dequantized}.11:} FDH Orientation Input}}}
\put(214.029,0){\makebox(0,0)[cc]
{{{\bf Figure \ref{Dequantized}.12:} FDH Orientation Output.}}}
\put(159,80.5){\circle*{1.5}}
\put(17,80.75){\circle*{1.5}}
\put(158.5,51){\circle*{1.5}}
\put(16.75,51.5){\circle*{1.5}}
\put(159.25,111){\circle*{1.5}}
\put(17.25,111.25){\circle*{1.5}}
\put(172.25,80.5){\circle*{1.5}}
\put(30.25,80.75){\circle*{1.5}}
\put(171.75,51){\circle*{1.5}}
\put(30,51.5){\circle*{1.5}}
\put(172.5,111){\circle*{1.5}}
\put(30.5,111.25){\circle*{1.5}}
\put(57.25,80.25){\circle*{1.5}}
\put(198.75,50.5){\circle*{1.5}}
\put(57,51.354){\circle*{1.5}}
\put(199.5,110.5){\circle*{1.5}}
\put(57.5,110.75){\circle*{1.5}}
\put(225.75,80.5){\circle*{1.5}}
\put(83.75,80.75){\circle*{1.5}}
\put(225.25,51){\circle*{1.5}}
\put(83.5,51.5){\circle*{1.5}}
\put(226,111){\circle*{1.5}}
\put(84,111.25){\circle*{1.5}}
\put(158.75,66.75){\circle*{1.5}}
\put(16.75,67){\circle*{1.5}}
\put(158.25,37.25){\circle*{1.5}}
\put(16.854,38.28){\circle*{1.5}}
\put(158.75,24.25){\circle*{1.5}}
\put(17,24.75){\circle*{1.5}}
\put(159,97.25){\circle*{1.5}}
\put(17,97.5){\circle*{1.5}}
\put(172,66.75){\circle*{1.5}}
\put(30,67){\circle*{1.5}}
\put(171.116,37.737){\circle*{1.5}}
\put(29.75,38.457){\circle*{1.5}}
\put(172,24.25){\circle*{1.5}}
\put(30.25,24.75){\circle*{1.5}}
\put(172.25,97.25){\circle*{1.5}}
\put(30.25,97.5){\circle*{1.5}}
\put(199,66.25){\circle*{1.5}}
\put(57,66.5){\circle*{1.5}}
\put(198.5,36.75){\circle*{1.5}}
\put(56.75,37.25){\circle*{1.5}}
\put(199,23.75){\circle*{1.5}}
\put(57.25,24.25){\circle*{1.5}}
\put(199.25,96.75){\circle*{1.5}}
\put(198.75,80){\circle*{1.5}}
\put(57.25,97){\circle*{1.5}}
\put(225.5,66.75){\circle*{1.5}}
\put(83.5,66.108){\circle*{1.5}}
\put(225,37.25){\circle*{1.5}}
\put(83.25,37.75){\circle*{1.5}}
\put(225.5,24.25){\circle*{1.5}}
\put(83.75,24.75){\circle*{1.5}}
\put(225.75,97.25){\circle*{1.5}}
\put(83.75,97.5){\circle*{1.5}}
\put(172.25,80.25){\circle*{1.5}}
\put(30.25,80.5){\circle*{1.5}}
\put(171.75,50.75){\circle*{1.5}}
\put(30,51.25){\circle*{1.5}}
\put(172.5,110.75){\circle*{1.5}}
\put(30.5,111){\circle*{1.5}}
\put(185.5,80.25){\circle*{1.5}}
\put(43.5,81.25){\circle*{1.5}}
\put(185,50.75){\circle*{1.5}}
\put(43.25,51.25){\circle*{1.5}}
\put(185.75,110.75){\circle*{1.5}}
\put(43.75,111){\circle*{1.5}}
\put(212.5,79.75){\circle*{1.5}}
\put(70.854,80.707){\circle*{1.5}}
\put(212,50.25){\circle*{1.5}}
\put(70.427,51.28){\circle*{1.5}}
\put(212.75,110.25){\circle*{1.5}}
\put(70.75,110.5){\circle*{1.5}}
\put(239,80.25){\circle*{1.5}}
\put(97,80.5){\circle*{1.5}}
\put(238.5,50.75){\circle*{1.5}}
\put(96.75,51.25){\circle*{1.5}}
\put(239.25,110.75){\circle*{1.5}}
\put(97.25,111){\circle*{1.5}}
\put(172,66.5){\circle*{1.5}}
\put(30,66.75){\circle*{1.5}}
\put(171.116,37.487){\circle*{1.5}}
\put(29.75,38.207){\circle*{1.5}}
\put(172,24){\circle*{1.5}}
\put(30.25,24.5){\circle*{1.5}}
\put(172.25,97){\circle*{1.5}}
\put(30.25,97.25){\circle*{1.5}}
\put(186,66.25){\circle*{1.5}}
\put(43.25,66.75){\circle*{1.5}}
\put(184.75,37){\circle*{1.5}}
\put(43,37.5){\circle*{1.5}}
\put(185.25,24){\circle*{1.5}}
\put(43.5,24.5){\circle*{1.5}}
\put(185.5,97){\circle*{1.5}}
\put(43.5,97.25){\circle*{1.5}}
\put(212.25,66){\circle*{1.5}}
\put(70.25,66.25){\circle*{1.5}}
\put(211.75,36.5){\circle*{1.5}}
\put(70,37){\circle*{1.5}}
\put(212.25,23.5){\circle*{1.5}}
\put(70.5,24){\circle*{1.5}}
\put(212.146,96.854){\circle*{1.5}}
\put(70.854,97.634){\circle*{1.5}}
\put(238.75,66.5){\circle*{1.5}}
\put(96.75,66.75){\circle*{1.5}}
\put(238.25,37){\circle*{1.5}}
\put(96.5,37.5){\circle*{1.5}}
\put(238.75,24){\circle*{1.5}}
\put(97,24.5){\circle*{1.5}}
\put(239,97){\circle*{1.5}}
\put(97,97.25){\circle*{1.5}}
\put(212.5,43.5){\makebox(0,0)[cc]{$z_{\bar k}$}}
\put(73.378,43.854){\makebox(0,0)[cc]{$z_{\bar k}$}}
\put(186.75,112){\line(0,1){0}}
\put(226,16){\line(0,1){0}}
\put(17,11.75){\line(0,1){26.5}}
\put(158.5,10.75){\line(0,1){26.5}}
\put(17,38.25){\line(1,0){12.75}}
\put(158.5,37.25){\line(1,0){12.75}}
\put(29.75,38.25){\line(0,1){13.75}}
\put(171.25,37.25){\line(0,1){13.75}}
\multiput(171.25,51)(8.91667,-.08333){3}{\line(1,0){8.91667}}
\put(56.75,81){\line(1,0){13.75}}
\put(198.25,80){\line(1,0){13.75}}
\put(70.5,81){\line(0,1){15.75}}
\multiput(70.5,97.5)(8.83333,-.08333){3}{\line(1,0){8.83333}}
\put(97,97.25){\line(0,1){0}}
\put(70.75,51.5){\vector(-1,0){41.25}}
\put(212,97.25){\line(1,0){27}}
\multiput(198.746,80.034)(.05893,-9.54594){3}{\line(0,-1){9.54594}}
\put(57.147,80.387){\line(0,-1){28.638}}
\put(212.004,50.689){\vector(0,1){46.492}}
\end{picture}

}
\end{center}
\medskip

\subsection{Numerical Example}

In this subsection, we consider a  numerical example where the number of inefficient firms is large enough to compare the scores. Specifically, we examine a graph translation homothetic structure in all cases except for the DEA and FDH cases. The input and output translation distance functions are identical in the translation homothetic case. This is not true for the convex and FDH cases, where only the output distance functions are reported.

\begin{center}{\scriptsize Table 1. Data Sample\\*{ \begin{tabular}
{|p{0.6cm}||p{0.9cm}|p{1.2cm}|p{1.2cm}|} \hline Firms & Input 1&
Input 2& Output \\ \hline \hline
1&1& 3&2\\ \hline 2 & 2& 2& 2\\
\hline 3
&2&1&2\\
\hline 4 &1&3&3\\
\hline
5 & 1& 4& 2\\
\hline
6 & 3& 2& 3\\
\hline 7 &4& 4& 5\\\hline
\end{tabular}}}
\end{center}

The values of the efficiency measures are listed for any value of $\alpha$
and are reported in the following table. Notably, it is clear that the technology is essentially characterized by units 4 and 7, which are efficient in every case.

\begin{center}{ \scriptsize Table 2. Efficiency scores under a VRS assumption.
{
\begin{tabular}
{|p{0.6cm}||p{1.2cm}|p{1.1cm}|p{1.1cm}| p{1.1cm}|p{1.1cm}|p{1.1 cm}| p{1.1cm}|p{1.1 cm}|p{1.1 cm}| p{1.2cm}|}
\hline Firms & $\alpha=-\infty$ &
$\alpha=-2$ &$\alpha=-1$ &$\alpha=-\frac{1}{2}$ & FDH& Convex& $\alpha=\frac{1}{2}$&$\alpha=1$& $\alpha=2$& $\alpha=+\infty$ 
  \\
\hline \hline
1		&1	&0.9729 	&0.9932		&1		&1		&1			&1				&1 				&0.9830	&1	\\
\hline
2 		&2	&1.6836 	&1.5475		&1.4445	&0		&0.8889		&1.2282			&1.14383		&1	 	&1\\
\hline
3		&1	& 0.9729  	&0.9932		&1 		&0		&0			&1				&1 				& 1.0012&1\\
\hline
4		&0	&0			&0			&0		&0		&0			&0 				&0				&0		&0\\
\hline
5		&1	&0.9729 	&1.0986		&1		&1		&1			&1 				&1				&1		&1\\
\hline
6		&1	&0.9729 	&0.9932		&1		&0		&0.5		&1 				&1				& 1.0012&1 \\
\hline
 7		&0	&0			&0			&0		&0		&0 		&0				&0				&0		&0
\\
\hline
\end{tabular}}}
\end{center}

The analysis of the results clearly aligns with the theoretical statements established in the paper. Though there is no evidence that the limit of the distance function should coincide with the distance function computed with respect to the limit technology, the geometric deformation of the technology necessarily impacts the efficiency scores. In particular, one can easily see that when 
$\alpha=-2$, the scores are close to those obtained in the case of the tropical limit 
$\alpha=-\infty$. However, note that the exponential nature of the algebraic operations involved in the computation of efficiency scores may cause some numerical problems; this is why we limit our calculations to small values of 
$\alpha$. The situation is symmetrical for the case 
$\alpha \longrightarrow +\infty  $, where the efficiency scores converge to those obtained at the tropical limit $\alpha=+\infty$. For $\alpha=2$, the efficiency scores are very close to those obtained at the tropical limit. The main intuition is that the quantized production technologies converge very quickly to their tropical limit.
  One can also observe that in the tropical and FDH cases, the distance functions are integer-valued. This shows that these approaches may have some practical issues when analyzing situations where the data are discrete (this is a common situation when the data are of a qualitative nature). In our example, the FDH case exhibits only two inefficient firms, making individual comparison more difficult. This is again more problematic in the input-oriented case, where an eyeball examination shows that  the scores are all identical to 0.

\section{Conclusion}\label{conclusion} 

Maslov's dequantization principle has been applied to production technology sets. In particular, Kolm-Pollack's generalized mean has been shown to provide a relevant semi-ring endowed with a proper algebraic structure, which yields either Min-plus or Max-plus production models in the neighborhood of infinity, i.e., tropical technologies obtained with the Painlevé-Kuratowski limit. These classes of technology sets are well-suited to design a graph translation homothetic structure compatible with output and input distance functions, which measure the degree of technical efficiency of a given firm or group of firms.

Discrete technology sets have been introduced to measure efficiency with ordinal data defined on an integer scale. The resulting input and output distance functions are shown to be discrete as well, and these are shown to be consistent with the well-known FDH technology

\end{document}